\documentclass[12pt,twoside]{article}
\usepackage[T1]{fontenc}
\usepackage{textcomp}
\usepackage{indentfirst}
\usepackage{bm}
\usepackage{amsmath}
\usepackage{mathrsfs}
\usepackage{amssymb}
\usepackage{fancyhdr}
\usepackage{latexsym}
\usepackage{bbding}
\usepackage{mathrsfs}
\usepackage{wasysym}
\usepackage{cite}
\usepackage{multicol,graphics}
\usepackage{xcolor}
\usepackage{enumitem}
\usepackage{lineno,hyperref}
\RequirePackage{times}

\def\XXint#1#2#3{{\setbox0=\hbox{$#1{#2#3}{\int}$ }
		\vcenter{\hbox{$#2#3$ }}\kern-.6\wd0}}

\newtheorem{theorem}{Theorem}[section]
\newtheorem{lemma}[theorem]{Lemma}
\newtheorem{corollary}[theorem]{Corollary}

\newtheorem{remark}[theorem]{Remark}

\numberwithin{equation}{section}

\allowdisplaybreaks[4]
\makeatletter
\begin{document}
	
\title{Liouville Type Theorem for the Fractional MHD and Hall-MHD equations in $\mathbb{R}^{3}$\footnote{This work was supported by National Natural Science Foundation of China (Grant No. 12271470).}}
\author{Weihua Wang\footnote{Corresponding author: wangvh@163.com (W. Wang).} ~\footnote{Authors are listed alphabetically by surname then given name. Authors equally share the first authorship.}, Zhenyuan Liu\footnote{mx120230298@stu.yzu.edu.cn}\\
[0.2cm] {\small School of Mathematical Sciences, Yangzhou University,}\\
[0.2cm] {\small  Yangzhou, Jiangsu 225002,  China}}
\date{\today}
\maketitle
	\begin{abstract}
\indent In this paper, we are mainly concerned with the Liouville type problem for the stationary fractional magnetohydrodynamics(MHD) and stationary fractional Hall-MHD equations. In addition, we present the results of the Navier-Stokes equation as a byproduct. The key point is to use the Caffarelli-Sivestre extension to overcome the difficulty caused by the non-local operator $(-\triangle)^{s}$ and combined with Yuan and Xiao's method (J. Math. Anal. Appl. 491 (2020) 124343).
	\end{abstract}
	\smallbreak
	
	\smallskip
	
	
\section{Introduction}
In this paper, we investigate Liouville-type problems for the following  stationary fractional magnetohydrodynamics (MHD) and Hall-MHD equations in $\mathbb{R}^{3}$
\begin{equation}\label{eq1.1}\tag{fMHD}
\begin{cases}
    (-\triangle)^{\alpha}u+(u \cdot \nabla)u-(b \cdot \nabla)b+\nabla p = 0,\\
    (-\triangle)^{\beta}b+(u \cdot \nabla)b-(b \cdot \nabla)u = 0,\\
     \nabla \cdot u = \nabla \cdot b = 0
\end{cases}
\end{equation}
and
\begin{equation}\label{eq1.3}\tag{fHall-MHD}
\begin{cases}
    (-\triangle)^{\alpha}u+(u \cdot \nabla)u-(b \cdot \nabla)b+\nabla p = 0, \\
    (-\triangle)^{\beta}b+(u \cdot \nabla)b-(b \cdot \nabla)u=\nabla \times((\nabla\times b)\times b),\\
     \nabla \cdot u =  \nabla \cdot b = 0,
\end{cases}
\end{equation}
where $u, ~b$ and $p$ denote the velocity, magnetic, and scalar pressure fields, respectively.  Hall-MHD equations are derived strictly from Euler-Maxwell equations or kinetic model in \cite{MPA2011}, which plays an important role in many physical problems, such as magnetic reconnection in space plasmas \cite{TF1991}, star formation \cite{SC2001}, and also neutron stars \cite{ML1960}. The Hall term $\nabla \times((\nabla\times b)\times b)$ is derived from the Ohm's law and describes deviation from charge neutrality between the electrons and the ions. And the fractional operator $(-\triangle)^{\alpha},~(-\triangle)^{\beta}$ is defined at the Fourier multiplier by the symbol $|\xi|^{2\alpha}$ and $|\xi|^{2\beta}$, respectively.\\
\indent When $b=0$ and $\alpha=1$, the system \eqref{eq1.1} reduces to the classical stationary  Navier-Stokes equations. Since the pioneering work of Leray \cite{Leary1933}, several authors have studied the question under which conditions the only smooth solution of the Navier-Stokes equations  and  related equations in $\mathbb{R}^{3}$ is trivial zero, which is called Liouville type theorem. For example, Galdi \cite[Chapter X, Remark
X.9.4 and Theorem X.9.5, p.729]{Gal2011} obtained  that if the smooth solution $u\in \textit{L}^\frac{9}{2}(\mathbb{R}^3)$, then $u =0$. After a few years, this result has been improved in different settings: D. Chae and J. Wolf gave a logarithmic improvement in \cite{DCha2016}; And recently, Chamorro-Jarr\'{i}n-Lemari\'{e}-Rieusset \cite{DO2018} proved that if the solution $u\in L^2_{loc}(\mathbb{R}^3)$ verifies $u\in L^p(\mathbb{R}^3)$ with $3\leq p \leq\frac{9}{2}$, then $u=0$. In addition, Yuan and Xiao \cite{YX2020} improved the result of \cite{DO2018} to extend $p$ to $2\leq p <\frac{18}{5}$.\\
\indent When $b=0$, the system \eqref{eq1.1} reduces to the stationary fractional Navier-Stokes equations, that is,
\begin{equation}\label{eq1.2}\tag{fNS}
\begin{cases}
(-\triangle)^{\alpha}u+u\cdot \nabla u+\nabla p=0,\\
     \nabla \cdot u = 0.
\end{cases}
\end{equation}
Yang \cite{YJ2022} extend the well-known $L^{\frac{9}{2}}(\mathbb{R}^3)$ result for the classical Navier-Stokes equations to the fractional Navier-Stokes equations. They proposed that the smooth solution $u\in L^{\frac{9}{2}}(\mathbb{R}^3)$ to $\eqref{eq1.2}$ with $\frac{5}{6}\leq \alpha <1$ is trivial zero.

When $\alpha=1$, the system \eqref{eq1.1} reduces to the classical MHD equations. Schulz \cite{Sch2019} obtained a Liouville theorem for the steady \eqref{eq1.1} without the finite Dirichlet integral, which is if the smooth solution $(u,b)\in L^p(\mathbb{R}^3)\cap BMO^{-1}(\mathbb{R}^3)$with $p\in(2,6]$, then $u =b =0$. Later, Yuan and Xiao \cite{YX2020} derived a theorem, under the condition that if the smooth solution $(u,b)\in L^p(\mathbb{R}^3)$with $2\leq p\leq\frac{9}{2}$, then $u =b =0$.

In this paper, inspired by Yang \cite{YJ2022}, we use the Caffarelli-Sivestre extension \cite[Section 2.3]{CaSi2007} to convert the non-local operator $(-\triangle)^\alpha$ and $(-\triangle)^\beta$ on $\mathbb{R}^3$ to the local operator $\bar{\triangle}$ on $\mathbb{R}_+^4$, and then extend the result of \cite{YX2020} to fractional order case. In addition, we also extend Yang's results \cite{YJ2022} from the fractional  Navier-Stokes equations to the fractional  MHD and the fractional  Hall-MHD. Furthermore, compared with Zeng's results in \cite{zy2025}, our results fill the gap for the parameter range of $2 \leq p \leq 3$.\\
For the steady fractional MHD equations, we have the following results:
\begin{theorem}\label{th1.1}
 Let $\frac{1}{2}\leq\alpha,~\beta<1$ and $\Lambda^\alpha u, \Lambda^\beta b$ in $L^2(\mathbb{R}^3)$.  Assume that  $(u,b)$ is  a smooth solution to the system $\eqref{eq1.1}$, then $u=b=0$ provided $u, ~b \in L^p(\mathbb{R}^3)$ with $2\leq p\leq\frac{9}{2}$.
\end{theorem}
\begin{remark}
  In the proof, we apply H\"{o}lder inequality to obtain $$\|u\|_{\mathrm{L}^3(B_{\frac{3}{2}R}\backslash B_{\frac{3}{4}R})}\lesssim R^{3-\frac{9}{p}}\|u\|_{\mathrm{L}^p(B_{\frac{3}{2}R}\backslash B_{\frac{3}{4}R})}$$ in Subsection \ref{subsubsec3.1.1}, which requires
$p\geq3$. However, given the additional condition $\Lambda^\alpha u\in L^2(\mathbb{R}^3)$, we further deduce that $p\leq\frac{6}{3-2\alpha}$, because of $\dot{H}^\alpha(\mathbb{R}^3)\hookrightarrow L^{\frac{6}{3-2\alpha}}(\mathbb{R}^3)$.  Consequently, we have $3\leq\frac{6}{3-2\alpha}$, i.e. $\alpha\geq\frac{1}{2}$. In addition, we apply the Caffarelli-Silvestre extension, which requires $\alpha <1$. Combining this with the previous condition, the final range is determined to be $\frac{1}{2}\leq\alpha<1$. The same reasoning applies to the case of $\frac{1}{2}\leq\beta<1$.
\end{remark}
\begin{remark}
  It should be noted that Zeng investigated Liouville-type results for the fractional MHD equations in the context of anisotropic Lebesgue spaces in \cite{zy2025}. When $p_{j,1}=p_{j,2}=p_{j,3}$, the range he obtained was $3\leq p\leq\frac{9}{2}$, which is contained in the results in this paper.
\end{remark}
As a byproduct, we also provide the corresponding results of the steady fractional  Navier-Stokes  equations.
\begin{corollary}\label{col1.2}
By removing the condition on the magnetic field $b$ in Theorem \ref{th1.1} while retaining the condition on the velocity field $u$, we can obtain the Liouville theorem for the steady fractional  Navier-Stokes  equations \eqref{eq1.2}.
\end{corollary}
\begin{remark}
Compared with the Yuan and Xiao's results\cite[Theorem 1.1.]{YX2020} $2\leq p\leq\frac{18}{5}$ is contained in $2\leq p\leq\frac{9}{2}$ provided $\alpha=1$ in Corollary \ref{col1.2}.  And Yang \cite{YJ2022} proved the validity of Corollary \ref{col1.2} provided $\frac{5}{6}\leq \alpha < 1$ and $p=\frac{9}{2}$.
\end{remark}

When $\alpha=\beta=1$, the system \eqref{eq1.3} reduces to the classical Hall-MHD equations. The Liouville-type problem of the Hall-MHD equations  is studied recently, as shown in works \cite{{cp2014},{cw2016},{zx2015}}. Chae et al. proved that if a smooth solution $(u,b)\in L^\infty\cap L^{\frac{9}{2}}(\mathbb{R}^3)$ with satisfying a finite Dirichlet integral condition $(\nabla u,\nabla b)\in L^2(\mathbb{R}^3)$, then $u=b=0$. Zhang et al.\cite{zx2015} generalized Cai's result\cite{cp2014} to hold under $(u,b)\in L^{\frac{9}{2}}(\mathbb{R}^3),~(\nabla u, \nabla b) \in L^2(\mathbb{R}^3)$. Later, Chae and Weng \cite{cw2016}  established  Liouville-type result  for any solution $(u,b)$ to either system MHD or Hall-MHD  with $u\in L^3(\mathbb{R}^3)$ and $(\nabla u,\nabla b)\in L^2(\mathbb{R}^3)$. We extend the results of the aforementioned paper to the fractional case and broaden the range of the integrability exponent $p$. More precisely, our main results are as follows:
\begin{theorem}\label{th1.5}
Let $\frac{1}{2}\leq\alpha<1$,   $\frac{5}{6}\leq\beta < 1$ and $\Lambda^\alpha u,~\Lambda^\beta b $ in $ L^2(\mathbb{R}^3)$. Suppose that $(u,b)$ is a smooth solution to the system $\eqref{eq1.3}$, then $u=b=0$ provided $(u,b)\in L^p(\mathbb{R}^3)$ with $2\leq p\leq\frac{9}{2}$.
\end{theorem}

The plan of this paper is organized as follow: In Section 2, we introduce some cut-off functions and Caffarelli-Sivestre extension \cite[Section 2.3]{CaSi2007}, moreover, we also give some lemmas.
 Section 3 and Section 4 are devoted to proving Theorem \ref{th1.1} and Theorem \ref{th1.5}, respectively.

\section{Preliminaries}\label{Sec2}
Let $\|(u,b)\|_{X}^p$ stand for $\|u\|_{X}^p+\|b\|_{X}^p$.
Throughout this paper,   $C$  denotes a finite inessential constant which may be different from line to line and $B_R(x_0)$ always represents the ball of radius $R$ centered at $x_0$, with the special case $B_R:=B_R(0)$ when $x_0=0$.

Let $\psi_R(x)\in C_c^\infty(\mathbb{R}^3)$ be a cut-off function as following
$$\psi_R(x)=
\begin{cases}
  1, & x\in B_{\frac{3}{4}R};\\
  0, & x\in B_{\frac{3}{2}R}^c.
\end{cases}$$
Then, it is easy to know $\|\nabla^s\psi_R(x)\|_{L^\infty(\mathbb{R}^3)}\leq\frac{C}{R^s}$.

In addition, we set $\chi_R(y)$ be a real  smooth cut-off function defined in $\mathbb{R}_{+}\cup\{0\}$ as following
\begin{eqnarray*}
  \chi_R(y)=\left\{
  \begin{array}{ll}
    1, & 0\leq y\leq\frac{3}{4}R; \\
    0, & \hbox{$y>\frac{3}{2}R$},
  \end{array}
\right.
\end{eqnarray*}
then $$|\chi_R'(y)|\leq MR^{-1}, ~|\chi_R''(y)|\leq MR^{-2}$$  for some  constant $M$ independent of $y \in \mathbb{R}_{+}\cup\{0\}$.

Compared to Laplacian operators $\triangle$, in order to overcome the technical difficulties brought by fractional order Laplacian operators $(-\triangle)^{s}$,  we need to introduce "Caffarelli-Sivestre" extension operators\cite[Section 2.3]{CaSi2007}. We use the notation $\bar{\nabla},~\bar{\triangle}$ for the differential operators defined on the upper half-space $\mathbb{R}^{4}_{+}$.

According to \cite[Section 2.3]{CaSi2007} and  \cite[Theorem 2.3]{CMLC2020}, for any $\bm{w}\in\dot{H}^{s}(\mathbb{R}^{3})$ with $s$ in $(0,1)$,  there is a unique "extension" $\bm{w^*}$ of $\bm{w}$ in the weighted space $H^{1}(\mathbb{R}^{4}_{+},y^{\lambda_s})$ with $\lambda_s = 1-2s$ which satisfies
\begin{equation}\label{Eq2.1}
\begin{cases}
    \bar{\triangle}_{\lambda_s}\bm{w^*}(x,y) := \bar{\triangle}\bm{w^*}+\frac{\lambda_s}{y}\partial_{y}\bm{w^*} = \frac{1}{y^{\lambda_s}}\overline{div}(y^{\lambda_s}\bar{\nabla}\bm{w^*}) = 0,\\
    \bm{w^*}(x,0)=\bm{w}(x),\\
     (-\triangle)^{s}\bm{w}(x) = -\bm{C}_{s}\lim\limits_{y\rightarrow 0}y^{\lambda_s}\partial_{y}\bm{w^*}(x,y),\\
     \int_{\mathbb{R}^{3}}|(-\triangle)^{\frac{s}{2}}\bm{w}|^2 dx =  \int_{\mathbb{R}^{3}}|\xi|^{2s}|\hat{\bm{w}}(\xi)|^2 d\xi = \bm{C}_{s}\int_{\mathbb{R}^{4}_{+}}y^{\lambda_s}|\bar{\nabla}\bm{w^*}|^2 dxdy.
\end{cases}
\end{equation}
 where constant $\bm{C}_{s}$ depends  only on $s$.

 In the integer-order case,  $\nabla w = 0$  indicates that $w$ is a constant. Therefore, we can naturally extend this result to fractional order. Indeed, the following conclusion can be readily derived from the proof of  \cite[Theorem 2.1]{YJ2022}:
\begin{lemma}\label{le2.3}
\begin{equation*}
  \int_{\mathbb{R}^{4}_{+}}y^{\lambda_s}|\bar{\nabla}\bm{w^*}|^2 dxdy=0
\end{equation*}
 implies that $w^*$ is almost constant everywhere.
\end{lemma}
\begin{remark}
  According to \eqref{Eq2.1}, we have $\|\Lambda^{s}w\|_{L^{2}(\mathbb{R}^{3})}^{2} =\bm{C}_{s}\int_{\mathbb{R}^{4}_{+}}y^{\lambda_s}|\bar{\nabla}\bm{w^*}|^2 dxdy$. Hence, $\Lambda^{s}w =0$ implies that $w$ is almost constant everywhere by Lemma \ref{le2.3}.
\end{remark}

\begin{lemma}[Lemma 2.2 in \cite{WX2018}]\label{le2.5}
  Let $\alpha \in (0, 1)$ and $u^∗$ be the $\alpha$-extension of $u \in L^p(\mathbb{R}^3)$ given by \eqref{Eq2.1}, it follows that
  \begin{equation}\label{2.2}
     \left(\int\int_{R^4_+}y^{1-2\alpha}|u^*|^{\frac{(5-2\alpha)p}{3}} dxdy\right)^{\frac{3}{(5-2\alpha)p}}\leq C\|u\|_{L^p(\mathbb{R}^3)}. \\
  \end{equation}
  In view of the Sobolev embedding $\dot{H}^\alpha(\mathbb{R}^3)\hookrightarrow L^{\frac{6}{3-2\alpha}}(\mathbb{R}^3)$, if we let $p=\frac{6}{3-2\alpha}$, we can obtain
  \begin{equation}\label{2.3}
      \left(\int\int_{R^4_+}y^{1-2\alpha}|u^*|^{\frac{2(5-2\alpha)}{3-2\alpha}} dxdy\right)^{\frac{3-2\alpha}{2(5-2\alpha)}}\leq C\|u\|_{\dot{H}^\alpha(\mathbb{R}^3)}. \\
  \end{equation}
\end{lemma}

\begin{lemma}[Gagliardo-Nirenberg inequality,  ~Lemma 2.2, ~\cite{Wu24} or \cite{Nirenberg11}]\label{le2.4}
		Let $1\leq q,s<\infty$ and $m\leq k$.
		Suppose that $j$ and $\vartheta$ satisfy $m\leq j\leq k$, $0\leq \vartheta\leq 1$ and define $p\in [1,+\infty]$ by
		
		\begin{equation*}
			j- \frac{3}{p}=\left(m-\frac{3}{s}\right)\vartheta+\left(k-\frac{3}{q}\right)(1-\vartheta).
		\end{equation*}
		Then the inequality holds:
		\begin{equation*}
			\|\nabla^{j} u\|_{L^{p}}\leq C	\|\nabla^{m}u\|_{L^{s}}^{\vartheta}\|\nabla^{k} u\|_{L^{q}}^{1-\vartheta}, ~ u\in W^{m,s}(\mathbb{T}^{3})\cap W^{k,q}(\mathbb{T}^{3}),
		\end{equation*}
		where constant $C\geq 0$. Here, when $p=\infty$, we require that $0<\vartheta<1$.
\end{lemma}
	\begin{lemma}\label{le2.6}
		Suppose $f,g\in L^2(\mathbb{R}^3)$, then the Parseval identity holds:
		$$\int_{\mathbb{R}^{3}}f(x)\overline{g(x)}~dx=\int_{\mathbb{R}^{3}}\hat{f(\xi)}\overline{\hat{g(\xi)}}~d\xi$$
		
		Furthermore, in this paper, we have $\Lambda^\alpha f,~\Lambda^{1-\alpha} g\in L^2(\mathbb{R}^3)$, thus we can obtain:
		\begin{equation*}
			\int_{\mathbb{R}^3} \Lambda f \cdot g ~dx = \int_{\mathbb{R}^3} \Lambda^\alpha f \cdot \Lambda^{1-\alpha} g ~dx
		\end{equation*}
	\end{lemma}

\begin{lemma}[\cite{VA2019}, Fractional Leibniz rule]\label{le2.2}
Given that $f$ and $g$ are two smooth functions, then we have the following estimate:
  $$\|(-\triangle)^s(fg)\|_{L^p}\leq C\|(-\triangle)^sf\|_{L^{p_0}}\|g\|_{L^{p_1}}+C\|f\|_{L^{q_0}}\|(-\triangle)^s g\|_{L^{q_1}},$$
  provided $\frac{1}{p}=\frac{1}{p_0}+\frac{1}{p_1}=\frac{1}{q_0}+\frac{1}{q_1}$, with $s>0$, $1 < p < +\infty$, and $1 < p_0, p_1, q_0, q_1 \leq
+\infty$.
\end{lemma}

\section{The proof of Theorem \ref{th1.1}}

First, multiplying $\eqref{eq1.1}_1$ by $u\psi_R(x)\chi(y)$ and $\eqref{eq1.1}_2$ by $b\psi_R(x)\chi(y)$,
 \begin{eqnarray*}
   0 &=& \langle (-\triangle)^{\alpha}u\,\text{,}\,
   u\psi_{R}\rangle +\langle(u \cdot \nabla)u\,\text{,}\,u\psi_{R}\rangle-\langle(b \cdot \nabla)b \,\text{,}\,u\psi_{R}\rangle+\langle\nabla P\,\text{,}\,u\psi_{R} \rangle\\
   & &+\langle (-\triangle)^{\beta}u\,\text{,}\,b\psi_{R}\rangle+\langle(u \cdot \nabla)b\,\text{,}\,b\psi_{R}\rangle-\langle(b \cdot \nabla)u \,\text{,}\,b\psi_{R}\rangle.
 \end{eqnarray*}
  Integrating by parts over $\mathbb{R}^{3}$ and taking into account $\eqref{eq1.1}_3$, we have
\begin{enumerate}[label=(\Roman*)]
  \item
$\langle(u \cdot \nabla) u\,\text{,}\,\psi_R u \rangle = \int_{\mathbb{R}^{3}} u_i \partial_i u_j \psi_R u_j \, dx = \frac{1}{2} \int_{\mathbb{R}^{3}} u_i \partial_i u_j^2 \psi_R \, dx = - \int_{\mathbb{R}^{3}} (u \cdot \nabla \psi_R) \frac{|u|^2}{2} \, dx$
  \item\label{2}
$\langle(b \cdot \nabla) b\,\text{,}\,\psi_R u \, \rangle= \int_{\mathbb{R}^{3}} b_i \partial_i b_j \psi_R u_j \, dx ,$
\item
$\langle\nabla p\,\text{,}\,\psi_Ru\rangle=\int_{\mathbb{R}^{3}}\partial_ip\cdot \psi_Ru_i~dx=-\int_{\mathbb{R}^{3}}p\cdot \partial_i(\psi_Ru_i)~dx=-\int_{\mathbb{R}^{3}}p\partial_i\psi_Ru_i~dx\\
=-\int_{\mathbb{R}^{3}}pu\cdot\nabla\psi_R~dx,$
\item
$\langle(u\cdot\nabla)b\,\text{,}\,\psi_Rb\rangle=\int_{\mathbb{R}^{3}}u_i\partial_ib_j\cdot\psi_Rb_j~dx=\frac{1}{2}\int_{\mathbb{R}^{3}}u_i\partial_ib_j^2\psi_R~dx=-\frac{1}{2}\int_{\mathbb{R}^{3}}u_ib_j^2\partial_i\psi_R~dx\\
=-\int_{\mathbb{R}^{3}}(u\cdot\nabla\psi_R)\frac{|b|^2}{2}~dx,$
\item \label{5}
$\langle(b\cdot\nabla)u\,\text{,}\,\psi_Rb\rangle=\int_{\mathbb{R}^{3}}b_i\partial_iu_j\cdot\psi_Rb_j~dx=-\int_{\mathbb{R}^{3}}b_iu_j\partial_i(\psi_Rb_j)~dx=-\int_{\mathbb{R}^{3}}b_iu_j\partial_i\psi_Rb_j~dx\\
-\int_{\mathbb{R}^{3}}b_iu_j\psi_R\partial_ib_j~dx.$
\end{enumerate}
We notice that there are the same items in \ref{2} and \ref{5}, so we add the two equations together to obtain
$$\ref{2}+ \ref{5}=-\int_{\mathbb{R}^{3}}b_iu_j\partial_i\psi_Rb_j~dx=-\int_{\mathbb{R}^{3}}(b\cdot\nabla\psi_R)(b\cdot u).$$
So we have
\begin{align}\label{eq3.1}
    \langle(-\triangle)^{\alpha}u\,\text{,}\,\psi_Ru\rangle+\langle(-\triangle)^{\beta}b\,\text{,}\,\psi_Rb\rangle=&\int_{\mathbb{R}^{3}}(u\cdot\nabla\psi_R)\left(\frac{|u|^2}{2}+\frac{|b|^2}{2}\right) \nonumber- \int_{\mathbb{R}^{3}}(b\cdot\nabla\psi_R)(b\cdot u)~dx\nonumber\\
    &+\int_{\mathbb{R}^{3}}pu\cdot\nabla\psi_R~dx.
\end{align}
We first derive an estimate of the pressure. Taking the divergence of $\eqref{eq1.1}_1$,  we can obtain
  $$\nabla\cdot((u\cdot\nabla)u)+\triangle p =\nabla\cdot((b\cdot\nabla)b),$$
       then
\begin{eqnarray*}
        p&=&\sum_{i,j=1}^3-\frac{1}{\triangle}\partial_i\partial_j(u_iu_j-b_ib_j)\\
         &=&\sum_{i,j=1}^3\mathcal{R}_i\mathcal{R}_j(u_iu_j-b_ib_j),
\end{eqnarray*}
where $\mathcal{R}_i=\frac{\partial_i}{\sqrt{-\triangle}}$ denotes the i-th Riesz transform. By the boundedness of the operator $\mathcal{R}_i$ on Lebesgue space $\mathrm{L}^p$ $(1<p<+\infty)$, we have
\begin{equation}\label{eq3.5}
    \|p\|_{\mathrm{L}^{\frac{p}{2}}}\leq C\|u\|_{\mathrm{L}^p}^2+C\|b\|_{\mathrm{L}^p}^2,
\end{equation}
which implies that the pressure $p$ belongs to $\mathrm{L}^\frac{p}{2}(\mathbb{R}^3)$ provided that $u$ and $b$ belong to $\mathrm{L}^p(\mathbb{R}^3)$.

In order to estimate $(-\triangle)^\alpha u$ and $(-\triangle)^\beta b$, we multiply \eqref{Eq2.1}$_{1}$ by $C_\alpha y^{\lambda_\alpha}u^*\psi_R(x)\chi_R(y)$ and then integrate to obtain
  \begin{align}\label{eq3.2}
    0=&C_{\alpha}\int_{\mathbb{R}^{4}_+}\overline{div}(y^{\lambda_\alpha}\bar{\nabla}u^*)\cdot u^*(x,y)~\psi_R(x)\chi_R(y) dxdy\nonumber\\
    =&-C_{\alpha}\int_{\mathbb{R}^{3}}\lim\limits_{y\rightarrow 0}~y^{\lambda_\alpha}\partial_{y}u^*\cdot u^*(x,y)\psi_R(x)\chi_R(y) dx\nonumber\\
    &-C_{\alpha}\int_{\mathbb{R}^{4}_+}y^{\lambda_\alpha}|\bar{\nabla}u^*|^2~\psi_R(x)\chi_R(y)~ dxdy - C_{\alpha}\int_{\mathbb{R}^{4}_+}y^{\lambda_\alpha}\overline{\nabla}u^*\cdot u^*\bar{\nabla}(\psi_R(x)\chi_R(y)) dxdy\nonumber\\
    =&\langle(-\triangle)^\alpha u\,\text{,}\,\psi_Ru\rangle-C_{\alpha}\int_{\mathbb{R}^{4}_+}y^{\lambda_\alpha}|\bar{\nabla}u^*|^2~\psi_R(x)\chi_R(y) dxdy\nonumber\\
&  - C_{\alpha}\int_{\mathbb{R}^{4}_+}y^{\lambda_\alpha}\overline{\nabla}u^*\cdot u^*\bar{\nabla}(\psi_R(x)\chi_R(y)) dxdy.
\end{align}
Similar to Eq. \eqref{eq3.2}, we have
  \begin{align}\label{eq3.3}
    0=&C_{\beta}\int_{\mathbb{R}^{4}_+}\overline{div}(y^{\lambda_\beta}\bar{\nabla}b^*)\cdot b^*(x,y)~\psi_R(x)\chi_R(y) dxdy\nonumber\\
    =&-C_{\beta}\int_{\mathbb{R}^{3}}\lim\limits_{y\rightarrow 0}~y^{\lambda_\beta}\partial_{y}b^*\cdot b^*(x,y)\psi_R(x)\chi_R(y) dx\nonumber\\
    &-C_{\beta}\int_{\mathbb{R}^{4}_+}y^{\lambda_\beta}|\bar{\nabla}b^*|^2~\psi_R(x)\chi_R(y) dxdy - C_{\beta}\int_{\mathbb{R}^{4}_+}y^{\lambda_\beta}\bar{\nabla}b^*\cdot b^*\bar{\nabla}(\psi_R(x)\chi_R(y)) dxdy\nonumber\\
    =&\langle(-\triangle)^\beta b\,\text{,}\,\psi_Rb\rangle- C_{\beta}\int_{\mathbb{R}^{4}_+}y^{\lambda_\beta}|\bar{\nabla}b^*|^2~\psi_R(x)\chi_R(y) dxdy\nonumber\\
& - C_{\beta}\int_{\mathbb{R}^{4}_+}y^{\lambda_\beta}\bar{\nabla}b^*\cdot b^*\bar{\nabla}(\psi_R(x)\chi_R(y)) dxdy.
\end{align}
Combining the above two equations \eqref{eq3.2} - \eqref{eq3.3}, we can obtain that
  \begin{align*}
0=&\langle(-\triangle)^\alpha u\,\text{,}\,\psi_Ru\rangle+\langle(-\triangle)^\beta b\,\text{,}\,\psi_Rb\rangle\\
&-\bm{C}_{\alpha}\int_{\mathbb{R}^{4}_+}y^{\lambda_\alpha}|\bar{\nabla}u^*|^2~\psi_R(x)\chi_R(y) dxdy-\bm{C}_{\beta}\int_{\mathbb{R}^{4}_+}y^{\lambda_\beta}|\bar{\nabla}b^*|^2~\psi_R(x)\chi_R(y) dxdy\\
&- C_{\alpha}\int_{\mathbb{R}^{4}_+}y^{\lambda_\alpha}\overline{\nabla}u^*\cdot u^*\bar{\nabla}(\psi_R(x)\chi_R(y)) dxdy- C_{\beta}\int_{\mathbb{R}^{4}_+}y^{\lambda_\beta}\bar{\nabla}b^*\cdot b^*\bar{\nabla}(\psi_R(x)\chi_R(y)) dxdy.
\end{align*}
Since by substituting the equation above into Eq. \eqref{eq3.1}, we have
  \begin{align}\label{eq3.4}
&2\bm{C}_{\alpha}\int_{\mathbb{R}^{4}_+}y^{\lambda_\alpha}|\bar{\nabla}u^*|^2~\psi_R(x)\chi_R(y)~dxdy + 2\bm{C}_{\beta}y^{\lambda_\beta}|\bar{\nabla}b^*|^2~\psi_R(x)\chi_R(y)~dxdy\nonumber\\
=&\int_{\mathbb{R}^{3}}(u\cdot\nabla\psi_{R})|u|^2~dx + \int_{\mathbb{R}^{3}}(u\cdot\nabla\psi_{R})|b|^2~dx+2\int_{\mathbb{R}^{3}}(b\cdot u)(b\cdot \nabla\psi_R(x)) dx\nonumber\\
&+2\int_{\mathbb{R}^{3}}pu\cdot\nabla\psi_R~dx - C_{\alpha}\int_{\mathbb{R}^{4}_+}y^{\lambda_\alpha}\overline{\nabla}u^*\cdot u^*\bar{\nabla}(\psi_R(x)\chi_R(y)) dxdy\nonumber\\
&- C_{\beta}\int_{\mathbb{R}^{4}_+}y^{\lambda_\beta}\bar{\nabla}b^*\cdot b^*\bar{\nabla}(\psi_R(x)\chi_R(y)) dxdy\nonumber\\
:=& \sum_{i=1}^{6} I_{i},
\end{align}
Next, we will divide the proof of Theorem\ref{th1.1} into two parts.

\subsection{$\frac{3}{4}\leq\alpha, ~\beta<1$}

\subsubsection{Case $3\leq p \leq \frac{9}{2}$}\label{subsubsec3.1.1}
For the first term in Eq. \eqref{eq3.4}, we have
\begin{equation}\label{eq3.18}
 I_1=\int_{\mathbb{R}^{3}}(u\cdot\nabla\psi_{R})|u|^2~dx\lesssim \frac{1}{R}\|u\|_{\mathrm{L}^3(B_{\frac{3}{2}R}\backslash B_{\frac{3}{4}R})}^3\overset{(p\geq3)}{\lesssim} R^{2-\frac{9}{p}}\|u\|_{\mathrm{L}^p(B_{2R}\backslash B_{\frac{R}{2}})}^3.
\end{equation}
Similarly, we can obtain estimates for $I_2$ and $I_3$:
\begin{align}\label{eq3.19}
    I_2=&\int_{\mathbb{R}^{3}}(u\cdot\nabla\psi_{R})|b|^2~dx\nonumber\\
    \lesssim&\frac{1}{R}\|u\|_{\mathrm{L}^3(B_{\frac{3}{2}R}\backslash B_{\frac{3}{4}R})}\||b|^2\|_{\mathrm{L}^\frac{3}{2}(B_{\frac{3}{2}R}\backslash B_{\frac{3}{4}R})}\hbox{\qquad(By H\"{o}lder's inequality)}\nonumber\\
    \lesssim&\frac{1}{R}\|u\|_{\mathrm{L}^3(B_{2R}\backslash B_{\frac{R}{2}})}\|b\|_{\mathrm{L}^3(B_{2R}\backslash B_{\frac{R}{2}})}^2\nonumber\\
\lesssim&\frac{1}{R}\|(u,b)\|_{\mathrm{L}^3(B_{2R}\backslash B_{\frac{R}{2}})}^3\nonumber\\
\overset{(p\geq3)}{\lesssim}& R^{2-\frac{9}{p}}\|(u,b)\|_{\mathrm{L}^p(B_{2R}\backslash B_{\frac{R}{2}})}^3
\end{align}
and
\begin{align}\label{eq3.20}
I_3=&2\int_{\mathbb{R}^{3}}(b\cdot u)(b\cdot \nabla\psi_R(x)) dx\lesssim \int_{\mathbb{R}^{3}}|u||b|^2|\nabla\psi_R|~dx \nonumber\\
\lesssim&\frac{1}{R}\|(u,b)\|_{\mathrm{L}^3(B_{2R}\backslash B_{\frac{R}{2}})}^3\nonumber\\
\overset{(p\geq3)}{\lesssim}& R^{2-\frac{9}{p}}\|(u,b)\|_{\mathrm{L}^p(B_{2R}\backslash B_{\frac{R}{2}})}^3.
\end{align}
Next, we will use Eq. \eqref{eq3.5} to estimate $I_4$:
    \begin{align}\label{eq3.6}
        I_4=& 2\int_{\mathbb{R}^{3}}pu\cdot\nabla\psi_R~dx\nonumber\\
       \lesssim&\frac{1}{R}\|u\|_{\mathrm{L}^3(B_{\frac{3}{2}R}\backslash B_{\frac{3}{4}R})}\|p\|_{\mathrm{L}^{\frac{3}{2}}(B_{\frac{3}{2}R}\backslash B_{\frac{3}{4}R})}\hbox{\qquad(By H\"{o}lder's inequality)}\nonumber\\
        \lesssim&\frac{1}{R}\|u\|_{\mathrm{L}^3(B_{2R}\backslash B_{\frac{R}{2}})}\|(u,b)\|_{\mathrm{L}^3(\mathbb{R}^3)}^2\hbox{\qquad(By \eqref{eq3.5})}\nonumber\\
\overset{(p\geq3)}{\lesssim}& R^{2-\frac{9}{p}}\|u\|_{\mathrm{L}^p(B_{2R}\backslash B_{\frac{R}{2}})}\|(u,b)\|_{\mathrm{L}^p(\mathbb{R}^3)}^2.\hbox{\qquad(By H\"{o}lder's inequality)}
    \end{align}

For $I_5$, we notice that
\begin{align}\label{3.11}
I_5 =& - C_{\alpha}\int_{\mathbb{R}^{4}_+}y^{\lambda_\alpha}\overline{\nabla}u^*\cdot u^*\bar{\nabla}(\psi_R(x)\chi_R(y)) dxdy\nonumber\\
  \lesssim& \int_{\mathbb{R}^{4}_+}y^{\lambda_\alpha}\overline{\nabla}u^*\cdot u^*\nabla\psi_R(x)\chi_R(y) dxdy + \int_{\mathbb{R}^{4}_+}y^{\lambda_\alpha}\overline{\nabla}u^*\cdot u^*\psi_R(x)\bar{\nabla}\chi_R(y) dxdy\nonumber\\
    \lesssim&\frac{1}{R}\int_{0}^{\frac{3}{2}R}\int_{B_{\frac{3}{2}R}\backslash B_{\frac{3}{4}R}}y^{\lambda_\alpha}\overline{\nabla}u^*\cdot u^*~dxdy+ \frac{1}{R}\int_{\frac{3}{4}R}^{\frac{3}{2}R}\int_{B_{\frac{3}{2}R}}y^{\lambda_\alpha}\overline{\nabla}u^*\cdot u^*~dxdy.
 \end{align}
 Next, applying H\"{o}lder's inequality and Lemma \ref{le2.5}, we obtain
\begin{align}\label{3.10}
  I_5 \lesssim& \frac{1}{R}\left(\int_{0}^{\frac{3}{2}R}\int_{B_{\frac{3}{2}R}\backslash B_{\frac{3}{4}R}}y^{\lambda_\alpha}|\overline{\nabla}u^*|^2~dxdy\right)^\frac{1}{2}\left(\int_{0}^{\frac{3}{2}R}\int_{B_{\frac{3}{2}R}\backslash B_{\frac{3}{4}R}}y^{\lambda_\alpha}|u^*|^{\frac{2(5-2\alpha)}{3-2\alpha}}~dxdy\right)^\frac{3-2\alpha}{2(5-2\alpha)}\nonumber\\
  &\left(\int_{0}^{\frac{3}{2}R}\int_{B_{\frac{3}{2}R}\backslash B_{\frac{3}{4}R}}y^{\lambda_\alpha}~dxdy\right)^\frac{1}{5-2\alpha}\nonumber\\
  &+\frac{1}{R}\left(\int_{\frac{3}{4}R}^{\frac{3}{2}R}\int_{B_{\frac{3}{2}R}}y^{\lambda_\alpha}|\overline{\nabla}u^*|^2~dxdy\right)^\frac{1}{2}\left(\int_{\frac{3}{4}R}^{\frac{3}{2}R}\int_{B_{\frac{3}{2}R}}y^{\lambda_\alpha}|u^*|^{\frac{2(5-2\alpha)}{3-2\alpha}}~dxdy\right)^\frac{3-2\alpha}{2(5-2\alpha)}\nonumber\\
  &\left(\int_{\frac{3}{4}R}^{\frac{3}{2}R}\int_{B_{\frac{3}{2}R}}y^{\lambda_\alpha}~dxdy\right)^\frac{1}{5-2\alpha}\hbox{\qquad(By H\"{o}lder's inequality)}\nonumber\\
\lesssim &\|u\|_{\dot{H}^\alpha(\mathbb{R}^3)}\left(\int_{0}^{\frac{3}{2}R}\int_{B_{\frac{3}{2}R}\backslash B_{\frac{3}{4}R}}y^{\lambda_\alpha}|\overline{\nabla}u^*|^2~dxdy\right)^\frac{1}{2}\nonumber\\
  &+\|u\|_{\dot{H}^\alpha(\mathbb{R}^3)}\left(\int_{\frac{3}{4}R}^{\frac{3}{2}R}\int_{B_{\frac{3}{2}R}}y^{\lambda_\alpha}|\overline{\nabla}u^*|^2~dxdy\right)^\frac{1}{2}\hbox{\qquad(By Lemma \ref{le2.5})}
\end{align}
 In addition, since $\Lambda^\alpha u,~\Lambda^\beta b\in L^2(\mathbb{R}^3)$, we have
\begin{align}\label{3.19}
&\int_{0}^{\frac{3}{2}R}\int_{B_{\frac{3}{2}R}\backslash B_{\frac{3}{4}R}}y^{\lambda_\alpha}|\bar{\nabla}u^*|^2~dxdy\rightarrow 0,~~\int_{\frac{3}{4}R}^{\frac{3}{2}R}\int_{B_{\frac{3}{2}R}}y^{\lambda_\alpha}|\bar{\nabla}u^*|^2 dxdy\rightarrow 0\nonumber\\
\text{ and }~~~~~& \nonumber\\
&\int_{0}^{\frac{3}{2}R}\int_{B_{\frac{3}{2}R}\backslash B_{\frac{3}{4}R}}y^{\lambda_\beta}|\bar{\nabla}b^*|^2~dxdy\rightarrow0,
~~\int_{\frac{3}{4}R}^{\frac{3}{2}R}\int_{B_{\frac{3}{2}R}}y^{\lambda_\beta}|\bar{\nabla}b^*|^2 dxdy\rightarrow0,
\end{align}
as $R\rightarrow\infty$. Therefore, we can obtain $I_5\rightarrow 0$ as $R\rightarrow\infty$.

Finally, we only need to repeat the procedure in \eqref{3.11}~-~\eqref{3.10} by replacing $u$ with $b$ and $\alpha$ with $\beta$, then we can have the estimation of $I_6$
\begin{equation}\label{3.16}
  I_6 \lesssim \|b\|_{\dot{H}^\beta(\mathbb{R}^3)}\left(\int_{0}^{\frac{3}{2}R}\int_{B_{\frac{3}{2}R}\backslash B_{\frac{3}{4}R}}y^{\lambda_\beta}|\overline{\nabla}b^*|^2~dxdy\right)^\frac{1}{2}
  +\|b\|_{\dot{H}^\beta(\mathbb{R}^3)}\left(\int_{\frac{3}{4}R}^{\frac{3}{2}R}\int_{B_{\frac{3}{2}R}}y^{\lambda_\beta}|\overline{\nabla}b^*|^2~dxdy\right)^\frac{1}{2}.
\end{equation}
Similarly, it is obvious that $I_6\rightarrow 0$ as $R\rightarrow\infty$.

So, by summing \eqref{eq3.18}, \eqref{eq3.19}, \eqref{eq3.20}, \eqref{eq3.6}, \eqref{3.10} and \eqref{3.16}, we conclude that
  \begin{align}\label{eq3.13}
&2\bm{C}_{\alpha}\int_{\mathbb{R}^{4}_+}y^{\lambda_\alpha}|\bar{\nabla}u^*|^2~\psi_R(x)\chi_R(y)~dxdy + 2\bm{C}_{\beta}y^{\lambda_\beta}|\bar{\nabla}b^*|^2~\psi_R(x)\chi_R(y)~dxdy\nonumber\\
\leq&\sum_{i=1}^{6}I_i\nonumber\\
\lesssim&\underbrace{R^{2-\frac{9}{p}}\|(u,b)\|_{\mathrm{L}^p(B_{2R}\backslash B_{\frac{R}{2}})}^3}_{I_{1} + I_{2}+ I_{3} \lesssim}+\underbrace{R^{2-\frac{9}{p}}\|u\|_{\mathrm{L}^p(B_{2R}\backslash B_{\frac{R}{2}})}\|(u,b)\|_{\mathrm{L}^p(\mathbb{R}^3)}^2}_{I_{4}\lesssim} \nonumber\\
&+  \underbrace{\|u\|_{\dot{H}^\alpha(\mathbb{R}^3)}\left(\int_{0}^{\frac{3}{2}R}\int_{B_{\frac{3}{2}R}\backslash B_{\frac{3}{4}R}}y^{\lambda_\alpha}|\overline{\nabla}u^*|^2~dxdy\right)^\frac{1}{2}+\|u\|_{\dot{H}^\alpha(\mathbb{R}^3)}\left(\int_{\frac{3}{4}R}^{\frac{3}{2}R}\int_{B_{\frac{3}{2}R}}y^{\lambda_\alpha}|\overline{\nabla}u^*|^2~dxdy\right)^\frac{1}{2}}_{I_{5}\lesssim}\nonumber\\
    &+ \underbrace{\|b\|_{\dot{H}^\beta(\mathbb{R}^3)}\left(\int_{0}^{\frac{3}{2}R}\int_{B_{\frac{3}{2}R}\backslash B_{\frac{3}{4}R}}y^{\lambda_\beta}|\overline{\nabla}b^*|^2~dxdy\right)^\frac{1}{2}
  +\|b\|_{\dot{H}^\beta(\mathbb{R}^3)}\left(\int_{\frac{3}{4}R}^{\frac{3}{2}R}\int_{B_{\frac{3}{2}R}}y^{\lambda_\beta}|\overline{\nabla}b^*|^2~dxdy\right)^\frac{1}{2}}_{I_{6}\lesssim}.
  \end{align}

Therefore, RHS of \eqref{eq3.13} converges to zero as $R\rightarrow\infty$ provided
\begin{equation*}
       2-\frac{9}{p}\leq 0 \Leftrightarrow p \leq \frac{9}{2},
\end{equation*}
that is, $3\leq p \leq \frac{9}{2}$.
Finally, we can obtain
$$2\bm{C}_{\alpha}\int_{\mathbb{R}^{4}_+}y^{\lambda_\alpha}|\bar{\nabla}u^*|^2~dxdy + 2\bm{C}_{\beta}\int_{\mathbb{R}^{4}_+}y^{\lambda_\beta}|\bar{\nabla}b^*|^2~dxdy\leq0,$$ i.e., $\bm{C}_{\alpha}\int_{\mathbb{R}^{4}_+}y^{\lambda_\alpha}|\bar{\nabla}u^*|^2~dxdy = 0$ and $\bm{C}_{\beta}\int_{\mathbb{R}^{4}_+}y^{\lambda_\beta}|\bar{\nabla}b^*|^2~dxdy = 0$,
which implies that $u$ and $b$ are almost constant everywhere according to Lemma \ref{le2.3} and \eqref{Eq2.1}$_{4}$. Since $(u,b)\in\mathrm{L}^p(\mathbb{R}^3)$ and $(u,b)$ is a smooth solution of \eqref{eq1.1}, we have $u=b=0$.

\subsubsection{Case $2\leq p \leq \min\{\frac{9}{2}, \frac{18(2\alpha-1)}{14\alpha-9}, \frac{18(2\beta-1)}{14\beta-9}\}$}

Firstly,  let $\phi_R(x)\in C_c^\infty(\mathbb{R}^3)$ be a non-negative cut-off function (for example, ref. \cite{YX2020}) as following:
\begin{eqnarray*}
  \phi_R(x)=\left\{
  \begin{array}{ll}
    1, & x\in B_{\frac{3}{2}R}\backslash B_{\frac{3}{4}R}, \\
    0, & x\in B_{2R}^c\cup B_{\frac{R}{2}}.
  \end{array}
\right.
\end{eqnarray*}
 And it is a simple matter to $$\|\nabla^s\phi_R\|_{L^{\infty}}\leq\frac{C}{R^s}.$$

 We notice that the estimation of $I_5$ and $I_6$ are the same as Eq. \eqref{3.10} and Eq. \eqref{3.16} in Subsection \ref{subsubsec3.1.1}, hence we only  need to estimate the nonlinear terms $I_1, ~I_2, ~I_3 \text{ and } I_4$.\\
In order to obtain the estimation of $I_{1}$, we have to estimate $\|u\|_{\mathrm{L}^4(B_{\frac{3}{2}R}\backslash B_{\frac{3}{4}R})}$ as follows
\begin{eqnarray}\label{3.18}
    & & \|u\|_{\mathrm{L}^4(B_{\frac{3}{2}R}\backslash B_{\frac{3}{4}R})}\leq \|u\phi_R\|_{\mathrm{L}^4(\mathbb{R}^3)}\nonumber\\
    &\overset{(\alpha\ge\frac{3}{4})}{\leq}& C\|u\phi_R\|_{\mathrm{L}^2(\mathbb{R}^3)}^{1-\frac{3}{4\alpha}}\|\Lambda^\alpha(u\phi_R)\|_{\mathrm{L}^2(\mathbb{R}^3)}^{\frac{3}{4\alpha}}\hbox{\quad(by Gagliardo-Nirenberg inequality)}\nonumber\\
    &\overset{\hbox{(Lemma \ref{le2.2})}}{\leq}&C\|u\phi_R\|_{\mathrm{L}^2(\mathbb{R}^3)}^{1-\frac{3}{4\alpha}}\left(\|\Lambda^\alpha u\|_{\mathrm{L}^2(B_{2R}\backslash B_{\frac{R}{2}})}^{\frac{3}{4\alpha}}\|\phi_R\|_{\mathrm{L}^\infty(\mathbb{R}^3)}^{\frac{3}{4\alpha}}+\|u\|_{\mathrm{L}^2(B_{2R}\backslash B_{\frac{R}{2}})}^{\frac{3}{4\alpha}}
    \|\Lambda^\alpha\phi_R\|_{\mathrm{L}^\infty(\mathbb{R}^3)}^{\frac{3}{4\alpha}}\right)\nonumber\\
    &\leq&C\|u\|_{\mathrm{L}^2(B_{2R}\backslash B_{\frac{R}{2}})}^{1-\frac{3}{4\alpha}}\left(\|\Lambda^\alpha u\|_{\mathrm{L}^2(B_{2R}\backslash B_{\frac{R}{2}})}^{\frac{3}{4\alpha}}+\frac{1}{R^\frac{3}{4}}\|u\|_{\mathrm{L}^2(B_{2R}\backslash B_{\frac{R}{2}})}^{\frac{3}{4\alpha}}\right)\nonumber\\
    &\leq&C\|u\|_{\mathrm{L}^2(B_{2R}\backslash B_{\frac{R}{2}})}^{1-\frac{3}{4\alpha}}\|\Lambda^\alpha u\|_{\mathrm{L}^2(B_{2R}\backslash B_{\frac{R}{2}})}^{\frac{3}{4\alpha}}+\frac{C}{R^\frac{3}{4}}\|u\|_{\mathrm{L}^2(B_{2R}\backslash B_{\frac{R}{2}})}.
\end{eqnarray}
 Thus, with the help of the H\"{o}lder's inequality and  Eq. \eqref{3.18}, we have
    \begin{align}\label{eq3.22}
I_1=&\int_{\mathbb{R}^{3}}(u\cdot\nabla\psi_{R})|u|^2~dx\leq\frac{C}{R}\|u\|_{\mathrm{L}^2(B_{\frac{3}{2}R}\backslash B_{\frac{3}{4}R})}\|u\|_{\mathrm{L}^4(B_{\frac{3}{2}R}\backslash B_{\frac{3}{4}R})}^2\nonumber\\
\leq&\frac{C}{R}\|u\|_{\mathrm{L}^2(B_{2R}\backslash B_{\frac{R}{2}})}\left(\|u\|_{\mathrm{L}^2(B_{2R}\backslash B_{\frac{R}{2}})}^{2-\frac{3}{2\alpha}}\|\Lambda^\alpha u\|_{\mathrm{L}^2(B_{2R}\backslash B_{\frac{R}{2}})}^{\frac{3}{2\alpha}}+\frac{C}{R^\frac{3}{2}}\|u\|_{\mathrm{L}^2(B_{2R}\backslash B_{\frac{R}{2}})}^2\right)\nonumber\\
\leq&\frac{C}{R}\|u\|_{\mathrm{L}^2(B_{2R}\backslash B_{\frac{R}{2}})}^{\frac{6\alpha-3}{2\alpha}}\|\Lambda^\alpha u\|_{\mathrm{L}^2(B_{2R}\backslash B_{\frac{R}{2}})}^{\frac{3}{2\alpha}}+\frac{C}{R^\frac{5}{2}}\|u\|_{\mathrm{L}^2(B_{2R}\backslash B_{\frac{R}{2}})}^3\nonumber\\
\overset{p\geq2}{\lesssim}&R^{\frac{9(2\alpha-1)(1-\frac{2}{p})-4\alpha}{4\alpha}}\|u\|_{\mathrm{L}^p(B_{2R}\backslash B_{\frac{R}{2}})}^{\frac{6\alpha-3}{2\alpha}}\|\Lambda^\alpha u\|_{\mathrm{L}^2(B_{2R}\backslash B_{\frac{R}{2}})}^{\frac{3}{2\alpha}}+ R^{2-\frac{9}{p}}\|u\|_{\mathrm{L}^p(B_{2R}\backslash B_{\frac{R}{2}})}^3.
    \end{align}
Under the conditions of $\beta\geq\frac{3}{4}$, we continue in the fashion to obtain
\begin{equation*}
        \|b\|_{\mathrm{L}^4(B_{\frac{3}{2}R}\backslash B_{\frac{3}{4}R})}\leq \|b\phi_R\|_{\mathrm{L}^4(\mathbb{R}^3)}\leq C\|b\|_{\mathrm{L}^2(B_{2R}\backslash B_{\frac{R}{2}})}^{1-\frac{3}{4\beta}}\|\Lambda^\beta b\|_{\mathrm{L}^2(B_{2R}\backslash B_{\frac{R}{2}})}^{\frac{3}{4\beta}}+\frac{C}{R^\frac{3}{2}}\|b\|_{\mathrm{L}^2(B_{2R}\backslash B_{\frac{R}{2}})}
   \end{equation*}
and
\begin{align}\label{eq3.23}
I_2=&\int_{\mathbb{R}^{3}}(u\cdot\nabla\psi_{R})|b|^2~dx\nonumber\\
 \leq&\frac{C}{R}\|u\|_{\mathrm{L}^2(B_{\frac{3}{2}R}\backslash B_{\frac{3}{4}R})}\|b\|_{\mathrm{L}^4(B_{\frac{3}{2}R}\backslash B_{\frac{3}{4}R})}^2\hbox{\qquad(by H\"{o}lder's inequality)}\nonumber\\
\leq&\frac{C}{R}\|u\|_{\mathrm{L}^2(B_{2R}\backslash B_{\frac{R}{2}})}\left(\|b\|_{\mathrm{L}^2(B_{2R}\backslash B_{\frac{R}{2}})}^{2-\frac{3}{2\beta}}\|\Lambda^\beta b\|_{\mathrm{L}^2(B_{2R}\backslash B_{\frac{R}{2}})}^{\frac{3}{2\beta}}+\frac{1}{R^\frac{3}{2}}\|b\|_{\mathrm{L}^2(B_{2R}\backslash B_{\frac{R}{2}})}^2\right)\nonumber\\
\leq&\frac{C}{R}\|(u,b)\|_{\mathrm{L}^2(B_{2R}\backslash B_{\frac{R}{2}})}^{\frac{6\beta-3}{2\beta}}\|\Lambda^\beta b\|_{\mathrm{L}^2(B_{2R}\backslash B_{\frac{R}{2}})}^{\frac{3}{2\beta}}+\frac{C}{R^\frac{5}{2}}\|(u,b)\|_{\mathrm{L}^2(B_{2R}\backslash B_{\frac{R}{2}})}^3\nonumber\\
\overset{p\geq2}{\lesssim}&R^{\frac{9(2\beta-1)(1-\frac{2}{p})-4\beta}{4\beta}}\|(u,b)\|_{\mathrm{L}^p(B_{2R}\backslash B_{\frac{R}{2}})}^{\frac{6\beta-3}{2\beta}}\|\Lambda^\beta b\|_{\mathrm{L}^2(B_{2R}\backslash B_{\frac{R}{2}})}^{\frac{3}{2\beta}}+ R^{2-\frac{9}{p}}\|(u,b)\|_{\mathrm{L}^p(B_{2R}\backslash B_{\frac{R}{2}})}^3.
\end{align}
Next, we use the estimation of $I_{2}$ to estimate $I_{3}$:
\begin{align}\label{eq3.15}
I_3=&2\int_{\mathbb{R}^{3}}(b\cdot u)(b\cdot \nabla\psi_R(x)) dx\leq\frac{C}{R}\|u\|_{\mathrm{L}^2(B_{\frac{3}{2}R}\backslash B_{\frac{3}{4}R})}\|b\|_{\mathrm{L}^4(B_{\frac{3}{2}R}\backslash B_{\frac{3}{4}R})}^2\hbox{\qquad(by H\"{o}lder's inequality)}\nonumber\\
&\hbox{(Using the estimation process of $I_{2}$ ......)}\nonumber\\
\overset{p\geq2}{\lesssim}&R^{\frac{9(2\beta-1)(1-\frac{2}{p})-4\beta}{4\beta}}\|(u,b)\|_{\mathrm{L}^p(B_{2R}\backslash B_{\frac{R}{2}})}^{\frac{6\beta-3}{2\beta}}\|\Lambda^\beta b\|_{\mathrm{L}^2(B_{2R}\backslash B_{\frac{R}{2}})}^{\frac{3}{2\beta}}+ R^{2-\frac{9}{p}}\|(u,b)\|_{\mathrm{L}^p(B_{2R}\backslash B_{\frac{R}{2}})}^3.
\end{align}
Then we need to combine the estimation of Eq. \eqref{eq3.22} and Eq. \eqref{eq3.23} to estimate $I_4$
\begin{align}\label{eq3.16}
    I_4=&2\int_{\mathbb{R}^{3}}pu\cdot\nabla\psi_R~dx\leq\frac{C}{R}\|u\|_{\mathrm{L}^2(B_{\frac{3}{2}R}\backslash B_{\frac{3}{4}R})}\|p\|_{\mathrm{L}^{2}(B_{\frac{3}{2}R}\backslash B_{\frac{3}{4}R})}\nonumber\\
    \leq&\frac{C}{R}\|u\|_{\mathrm{L}^2(B_{\frac{3}{2}R}\backslash B_{\frac{3}{4}R})}\left(\|u\|_{\mathrm{L}^4(\mathbb{R}^3)}^2+\|b\|_{\mathrm{L}^4(\mathbb{R}^3)}^2\right)\nonumber\\
    \leq&\frac{C}{R}\|u\|_{\mathrm{L}^2(B_{2R}\backslash B_{\frac{R}{2}})}\|u\|_{\mathrm{L}^2(\mathbb{R}^3)}^{\frac{4\alpha-3}{2\alpha}}\|\Lambda^\alpha u\|_{\mathrm{L}^2(\mathbb{R}^3)}^{\frac{3}{2\alpha}}+\frac{C}{R}\|u\|_{\mathrm{L}^2(B_{2R}\backslash B_{\frac{R}{2}})}\|b\|_{\mathrm{L}^2(\mathbb{R}^3)}^{\frac{4\beta-3}{2\beta}}\|\Lambda^\beta b\|_{\mathrm{L}^2(\mathbb{R}^3)}^{\frac{3}{2\beta}}\nonumber\\
    \overset{p\geq2}{\lesssim}&\underbrace{R^{\frac{9(2\alpha-1)(1-\frac{2}{p})-4\alpha}{4\alpha}}\|u\|_{\mathrm{L}^p(B_{2R}\backslash B_{\frac{R}{2}})}\|u\|_{\mathrm{L}^p(\mathbb{R}^3)}^{\frac{4\alpha-3}{2\alpha}}\|\Lambda^\alpha u\|_{\mathrm{L}^2(\mathbb{R}^3)}^{\frac{3}{2\alpha}}}_{I_{4,1}}\nonumber\\
    &+\underbrace{R^{\frac{9(2\beta-1)(1-\frac{2}{p})-4\beta}{4\beta}}\|u\|_{\mathrm{L}^p(B_{2R}\backslash B_{\frac{R}{2}})}\|b\|_{\mathrm{L}^p(\mathbb{R}^3)}^{\frac{4\beta-3}{2\beta}}\|\Lambda^\beta b\|_{\mathrm{L}^2(\mathbb{R}^3)}^{\frac{3}{2\beta}}}_{I_{4,2}}.
\end{align}
Then we combine $I_1-I_4$ \eqref{eq3.22}, \eqref{eq3.23}, \eqref{eq3.15}, \eqref{eq3.16} and $I_5$, $I_6$ \eqref{3.10}, \eqref{3.16} to obtain
\begin{align}\label{3.22}
&2\bm{C}_{\alpha}\int_{\mathbb{R}^{4}_+}y^{\lambda_\alpha}|\bar{\nabla}u^*|^2~\psi_R(x)\chi_R(y)~dxdy + 2\bm{C}_{\beta}y^{\lambda_\beta}|\bar{\nabla}b^*|^2~\psi_R(x)\chi_R(y)~dxdy\leq\sum_{i=1}^{6}I_i\nonumber\\
\lesssim&\underbrace{R^{\frac{9(2\alpha-1)(1-\frac{2}{p})-4\alpha}{4\alpha}}\|u\|_{\mathrm{L}^p(B_{2R}\backslash B_{\frac{R}{2}})}^{\frac{6\alpha-3}{2\alpha}}\|\Lambda^\alpha u\|_{\mathrm{L}^2(B_{2R}\backslash B_{\frac{R}{2}})}^{\frac{3}{2\alpha}}+ R^{2-\frac{9}{p}}\|u\|_{\mathrm{L}^p(B_{2R}\backslash B_{\frac{R}{2}})}^3}_{I_1\lesssim}\nonumber\\
&+\underbrace{R^{\frac{9(2\beta-1)(1-\frac{2}{p})-4\beta}{4\beta}}\|(u,b)\|_{\mathrm{L}^p(B_{2R}\backslash B_{\frac{R}{2}})}^{\frac{6\beta-3}{2\beta}}\|\Lambda^\beta b\|_{\mathrm{L}^2(B_{2R}\backslash B_{\frac{R}{2}})}^{\frac{3}{2\beta}}+ R^{2-\frac{9}{p}}\|(u,b)\|_{\mathrm{L}^p(B_{2R}\backslash B_{\frac{R}{2}})}^3}_{I_2+I_3\lesssim}\nonumber\\
&+I_{4,1}+I_{4,2}\nonumber\\
&+  \underbrace{\|u\|_{\dot{H}^\alpha(\mathbb{R}^3)}\left(\int_{0}^{\frac{3}{2}R}\int_{B_{\frac{3}{2}R}\backslash B_{\frac{3}{4}R}}y^{\lambda_\alpha}|\overline{\nabla}u^*|^2~dxdy\right)^\frac{1}{2}+\|u\|_{\dot{H}^\alpha(\mathbb{R}^3)}\left(\int_{\frac{3}{4}R}^{\frac{3}{2}R}\int_{B_{\frac{3}{2}R}}y^{\lambda_\alpha}|\overline{\nabla}u^*|^2~dxdy\right)^\frac{1}{2}}_{I_{5}\lesssim}\nonumber\\
    &+ \underbrace{\|b\|_{\dot{H}^\beta(\mathbb{R}^3)}\left(\int_{0}^{\frac{3}{2}R}\int_{B_{\frac{3}{2}R}\backslash B_{\frac{3}{4}R}}y^{\lambda_\beta}|\overline{\nabla}b^*|^2~dxdy\right)^\frac{1}{2}
  +\|b\|_{\dot{H}^\beta(\mathbb{R}^3)}\left(\int_{\frac{3}{4}R}^{\frac{3}{2}R}\int_{B_{\frac{3}{2}R}}y^{\lambda_\beta}|\overline{\nabla}b^*|^2~dxdy\right)^\frac{1}{2}}_{I_{6}\lesssim}\nonumber\\
\lesssim& R^{\frac{9(2\alpha-1)(1-\frac{2}{p})-4\alpha}{4\alpha}}\|u\|_{\mathrm{L}^p(B_{2R}\backslash B_{\frac{R}{2}})}^{3-\frac{3}{2\alpha}}\|\Lambda^\alpha u\|_{\mathrm{L}^2(B_{2R}\backslash B_{\frac{R}{2}})}^{\frac{3}{2\alpha}}\nonumber\\
&+ R^{\frac{9(2\beta-1)(1-\frac{2}{p})-4\beta}{4\beta}}\|(u,b)\|_{\mathrm{L}^p(B_{2R}\backslash B_{\frac{R}{2}})}^{\frac{6\beta-3}{2\beta}}\|\Lambda^\beta b\|_{\mathrm{L}^2(B_{2R}\backslash B_{\frac{R}{2}})}^{\frac{3}{2\beta}}+ R^{2-\frac{9}{p}}\|(u,b)\|_{\mathrm{L}^p(B_{2R}\backslash B_{\frac{R}{2}})}^3\nonumber\\
&+R^{\frac{9(2\alpha-1)(1-\frac{2}{p})-4\alpha}{4\alpha}}\|u\|_{\mathrm{L}^p(B_{2R}\backslash B_{\frac{R}{2}})}\|u\|_{\mathrm{L}^p(\mathbb{R}^3)}^{\frac{4\alpha-3}{2\alpha}}\|\Lambda^\alpha u\|_{\mathrm{L}^2(\mathbb{R}^3)}^{\frac{3}{2\alpha}}\nonumber\\
&+R^{\frac{9(2\beta-1)(1-\frac{2}{p})-4\beta}{4\beta}}\|u\|_{\mathrm{L}^p(B_{2R}\backslash B_{\frac{R}{2}})}\|b\|_{\mathrm{L}^p(\mathbb{R}^3)}^{\frac{4\beta-3}{2\beta}}\|\Lambda^\beta b\|_{\mathrm{L}^2(\mathbb{R}^3)}^{\frac{3}{2\beta}}\nonumber\\
&+\|u\|_{\dot{H}^\alpha(\mathbb{R}^3)}\left(\int_{0}^{\frac{3}{2}R}\int_{B_{\frac{3}{2}R}\backslash B_{\frac{3}{4}R}}y^{\lambda_\alpha}|\overline{\nabla}u^*|^2~dxdy\right)^\frac{1}{2}+\|u\|_{\dot{H}^\alpha(\mathbb{R}^3)}\left(\int_{\frac{3}{4}R}^{\frac{3}{2}R}\int_{B_{\frac{3}{2}R}}y^{\lambda_\alpha}|\overline{\nabla}u^*|^2~dxdy\right)^\frac{1}{2}\nonumber\\
    &+ \|b\|_{\dot{H}^\beta(\mathbb{R}^3)}\left(\int_{0}^{\frac{3}{2}R}\int_{B_{\frac{3}{2}R}\backslash B_{\frac{3}{4}R}}y^{\lambda_\beta}|\overline{\nabla}b^*|^2~dxdy\right)^\frac{1}{2}
  +\|b\|_{\dot{H}^\beta(\mathbb{R}^3)}\left(\int_{\frac{3}{4}R}^{\frac{3}{2}R}\int_{B_{\frac{3}{2}R}}y^{\lambda_\beta}|\overline{\nabla}b^*|^2~dxdy\right)^\frac{1}{2}\\
&\hbox{\quad(By H\"{o}lder's inequality)}\nonumber
  \end{align}
for $p\geq2$.\\
Considering $\Lambda^\alpha u\in L^2(\mathbb{R}^3)$ and \eqref{3.19}, RHS of Eq. \eqref{3.22} converges to zero as $R\rightarrow\infty$ when $p,~\alpha,~\beta$ satisfy
\begin{equation}\nonumber
  \left\{
     \begin{array}{ll}
     \frac{9(2\alpha-1)(1-\frac{2}{p})-4\alpha}{4\alpha} \leq 0  &  \hbox{$\Leftrightarrow p \leq \frac{18(2\alpha-1)}{14\alpha-9}$,}\\
     \frac{9(2\beta-1)(1-\frac{2}{p})-4\beta}{4\beta} \leq 0  &  \hbox{$\Leftrightarrow p \leq \frac{18(2\beta-1)}{14\beta-9}$,}\\
       2-\frac{9}{p}\leq 0 & \hbox{$\Leftrightarrow p \leq \frac{9}{2}$,} \\
     \end{array}
   \right.
\end{equation}
that is, $2\leq p \leq \min\{\frac{9}{2}, \frac{18(2\alpha-1)}{14\alpha-9}, \frac{18(2\beta-1)}{14\beta-9}\}$.

Therefore, it follows easily to obtain $\bm{C}_{\alpha}\int_{\mathbb{R}^{4}_+}y^{\lambda_\alpha}|\bar{\nabla}u^*|^2~dxdy = \bm{C}_{\beta}\int_{\mathbb{R}^{4}_+}y^{\lambda_\beta}|\bar{\nabla}b^*|^2~dxdy = 0$. According to the property stated in Lemma \ref{le2.3} and $\eqref{Eq2.1}_4$, since $(u,b)$ is a smooth solution of \eqref{eq1.1} in $L^p(\mathbb{R}^3)$, we have  $u=b\equiv0$.

It is easy to check that for $\alpha,\beta\geq\frac{3}{4}$, the inequalities $\frac{18(2\alpha-1)}{14\alpha-9} > 3$, $\frac{18(2\beta-1)}{14\beta-9}> 3$, $\frac{6}{3-2\alpha}>3$ and $\frac{6}{3-2\beta}> 3$ always hold true. Hence, by summarizing  Subsection 3.1.1 and Subsection 3.1.2, we can obtain $u=b=0$ provided $2\leq p\leq\frac{9}{2}$ with $\frac{3}{4}\leq\alpha,\beta<1$.\\
\subsection{$\frac{1}{2}\leq\alpha< \frac{3}{4}, ~\frac{1}{2}\leq\beta<1$}\label{subsec3.2}
 For the case of $3\leq p \leq \frac{9}{2}$, we only repeat the procedures in Subsection \ref{subsubsec3.1.1}. Therefore,  we have to proceed with the analysis for the case of $ 2\leq p\leq3 $.  Moreover, the estimations of $I_{5}$ and $I_{6}$  are the same as Subsection \ref{subsubsec3.1.1}.  The difference is that we need to modify the estimations of $I_1-I_4$.

Firstly, we apply the Gagliardo-Nirenberg inequality to estimate $\|u\|_{\mathrm{L}^3(B_{\frac{3}{2}R}\backslash B_{\frac{3}{4}R})}$:
\begin{align}\label{3.26}
 &\|u\|_{\mathrm{L}^3(B_{\frac{3}{2}R}\backslash B_{\frac{3}{4}R})}\leq \|u\phi_R\|_{\mathrm{L}^3(\mathbb{R}^3)}\nonumber\\
\leq& C\|u\phi_R\|_{\mathrm{L}^2(\mathbb{R}^3)}^{1-\frac{1}{2\alpha}}\|\Lambda^\alpha(u\phi_R)\|_{\mathrm{L}^2(\mathbb{R}^3)}^{\frac{1}{2\alpha}}\hbox{\quad (by Gagliardo-Nirenberg inequality)}\nonumber\\
    \leq&C\|u\|_{\mathrm{L}^2(B_{2R}\backslash B_{\frac{R}{2}})}^{1-\frac{1}{2\alpha}}\left(\|\Lambda^\alpha u\|_{\mathrm{L}^2(B_{2R}\backslash B_{\frac{R}{2}})}^{\frac{1}{2\alpha}}\|\phi_R\|_{\mathrm{L}^\infty(\mathbb{R}^3)}^{\frac{1}{2\alpha}}+
    \|u\|_{\mathrm{L}^2(B_{2R}\backslash B_{\frac{R}{2}})}^{\frac{1}{2\alpha}}\|\Lambda^\alpha\phi_R\|_{\mathrm{L}^\infty(\mathbb{R}^3)}^{\frac{1}{2\alpha}}\right)\nonumber\\
    \leq&C\|u\|_{\mathrm{L}^2(B_{2R}\backslash B_{\frac{R}{2}})}^{1-\frac{1}{2\alpha}}\left(\|\Lambda^\alpha u\|_{\mathrm{L}^2(B_{2R}\backslash B_{\frac{R}{2}})}^{\frac{1}{2\alpha}}+\frac{1}{R^\frac{1}{2}}\|u\|_{\mathrm{L}^2(B_{2R}\backslash B_{\frac{R}{2}})}^{\frac{1}{2\alpha}}\right)\nonumber\\
    \leq&C\|u\|_{\mathrm{L}^2(B_{2R}\backslash B_{\frac{R}{2}})}^{1-\frac{1}{2\alpha}}\|\Lambda^\alpha u\|_{\mathrm{L}^2(B_{2R}\backslash B_{\frac{R}{2}})}^{\frac{1}{2\alpha}}+\frac{C}{R^\frac{1}{2}}\|u\|_{\mathrm{L}^2(B_{2R}\backslash B_{\frac{R}{2}})}.
\end{align}

Similar considerations apply to $\|b\|_{\mathrm{L}^3(B_{\frac{3}{2}R}\backslash B_{\frac{3}{4}R})}$
\begin{equation}\label{3.31}
 \|b\|_{\mathrm{L}^3(B_{\frac{3}{2}R}\backslash B_{\frac{3}{4}R})}\leq C\|b\|_{\mathrm{L}^2(B_{2R}\backslash B_{\frac{R}{2}})}^{1-\frac{1}{2\beta}}\|\Lambda^\beta b\|_{\mathrm{L}^2(B_{2R}\backslash B_{\frac{R}{2}})}^{\frac{1}{2\beta}}+\frac{C}{R^\frac{1}{2}}\|b\|_{\mathrm{L}^2(B_{2R}\backslash B_{\frac{R}{2}})}.
\end{equation}

Then the estimation of $I_1$ is straightforward
\begin{align}\label{3.29}
I_1=&\int_{\mathbb{R}^{3}}(u\cdot\nabla\psi_{R})|u|^2~dx\leq\frac{C}{R}\|u\|_{\mathrm{L}^3(B_{\frac{3}{2}R}\backslash B_{\frac{3}{4}R})}^3\nonumber\\
\leq&\frac{C}{R}\|u\|_{\mathrm{L}^2(B_{2R}\backslash B_{\frac{R}{2}})}^{3-\frac{3}{2\alpha}}\|\Lambda^\alpha u\|_{\mathrm{L}^2(B_{2R}\backslash B_{\frac{R}{2}})}^{\frac{3}{2\alpha}}+\frac{C}{R^{\frac{5}{2}}}\|u\|_{\mathrm{L}^2(B_{2R}\backslash B_{\frac{R}{2}})}^3\nonumber\\
\overset{p\geq2}{\lesssim}&R^{\frac{3(1-\frac{2}{p})(6\alpha-3)}{4\alpha}-1}\|u\|_{\mathrm{L}^p(B_{2R}\backslash B_{\frac{R}{2}})}^{3-\frac{3}{2\alpha}}\|\Lambda^\alpha u\|_{\mathrm{L}^2(B_{2R}\backslash B_{\frac{R}{2}})}^{\frac{3}{2\alpha}}+ R^{2-\frac{9}{p}}\|u\|_{\mathrm{L}^p(B_{2R}\backslash B_{\frac{R}{2}})}^3,
\end{align}
According to \eqref{3.29}, it is necessary to choose
\begin{equation}\label{e1}
  \left\{
     \begin{array}{ll}
     \frac{3(1-\frac{2}{p})(6\alpha-3)}{4\alpha}-1\leq 0 , \\
     2-\frac{9}{p}\leq 0 . \\
     \end{array}
   \right.
\end{equation}
Hence, we will obtain $I_1 \rightarrow 0$ as $R \rightarrow \infty$.

Then we apply Eq. \eqref{3.26} and Eq. \eqref{3.31} to estimate $I_2$
\begin{align}\label{3.32}
  I_2=&\int_{\mathbb{R}^{3}}(u\cdot\nabla\psi_{R})|b|^2~dx\leq\frac{C}{R}\|u\|_{\mathrm{L}^3(B_{\frac{3}{2}R}\backslash B_{\frac{3}{4}R})} \|b\|_{\mathrm{L}^3(B_{\frac{3}{2}R}\backslash B_{\frac{3}{4}R})}^2 \nonumber\\
  \leq&\frac{C}{R}\left(\|u\|_{\mathrm{L}^2(B_{2R}\backslash B_{\frac{R}{2}})}^{1-\frac{1}{2\alpha}}\|\Lambda^\alpha u\|_{\mathrm{L}^2(B_{2R}\backslash B_{\frac{R}{2}})}^{\frac{1}{2\alpha}}+\frac{C}{R^\frac{1}{2}}\|u\|_{\mathrm{L}^2(B_{2R}\backslash B_{\frac{R}{2}})}\right)\nonumber\\
  &~~\qquad\left(\|b\|_{\mathrm{L}^2(B_{2R}\backslash B_{\frac{R}{2}})}^{2-\frac{1}{\beta}}\|\Lambda^\beta b\|_{\mathrm{L}^2(B_{2R}\backslash B_{\frac{R}{2}})}^{\frac{1}{\beta}}+\frac{C}{R}\|b\|_{\mathrm{L}^2(B_{2R}\backslash B_{\frac{R}{2}})}^2\right)\nonumber\\
  \leq&\frac{C}{R}\left(\|u\|_{\mathrm{L}^2(B_{2R}\backslash B_{\frac{R}{2}})}^{2-\frac{1}{\alpha}}\|\Lambda^\alpha u\|_{\mathrm{L}^2(B_{2R}\backslash B_{\frac{R}{2}})}^{\frac{1}{\alpha}}+\|b\|_{\mathrm{L}^2(B_{2R}\backslash B_{\frac{R}{2}})}^{4-\frac{2}{\beta}}\|\Lambda^\beta b\|_{\mathrm{L}^2(B_{2R}\backslash B_{\frac{R}{2}})}^{\frac{2}{\beta}}\right)\nonumber\\
  &+\frac{C}{R^2}\|u\|_{\mathrm{L}^2(B_{2R}\backslash B_{\frac{R}{2}})}^{2-\frac{1}{\alpha}}\|\Lambda^\alpha u\|_{\mathrm{L}^2(B_{2R}\backslash B_{\frac{R}{2}})}^{\frac{1}{\alpha}}+\frac{C}{R^2}\|b\|_{\mathrm{L}^2(B_{2R}\backslash B_{\frac{R}{2}})}^4\nonumber\\
  &+\frac{C}{R^{\frac{3}{2}}}\|(u,b)\|_{\mathrm{L}^2(B_{2R}\backslash B_{\frac{R}{2}})}^{3-\frac{1}{\beta}}\|\Lambda^\beta b\|_{\mathrm{L}^2(B_{2R}\backslash B_{\frac{R}{2}})}^{\frac{1}{\beta}}+\frac{C}{R^{\frac{5}{2}}}\|(u,b)\|_{\mathrm{L}^2(B_{2R}\backslash B_{\frac{R}{2}})}^{3}.\nonumber\\
  &\hbox{\qquad\qquad\qquad\qquad\qquad\qquad\qquad\qquad\qquad\qquad\qquad(by Young's inequality)}\nonumber\\
  \overset{(p\geq2)}{\lesssim}&R^{\frac{3(1-\frac{2}{p})(2\alpha-1)-2\alpha}{2\alpha}}\|u\|_{\mathrm{L}^p(B_{2R}\backslash B_{\frac{R}{2}})}^{2-\frac{1}{\alpha}}\|\Lambda^\alpha u\|_{\mathrm{L}^2(B_{2R}\backslash B_{\frac{R}{2}})}^{\frac{1}{\alpha}}\nonumber\\
&+ R^{\frac{3(1-\frac{2}{p})(2\beta-1)-\beta}{\beta}}\|b\|_{\mathrm{L}^p(B_{2R}\backslash B_{\frac{R}{2}})}^{4-\frac{2}{\beta}}\|\Lambda^\beta b\|_{\mathrm{L}^2(B_{2R}\backslash B_{\frac{R}{2}})}^{\frac{2}{\beta}}\nonumber\\
&+ R^{\frac{3(1-\frac{2}{p})(2\alpha-1)-4\alpha}{2\alpha}}\|u\|_{\mathrm{L}^p(B_{2R}\backslash B_{\frac{R}{2}})}^{2-\frac{1}{\alpha}}\|\Lambda^\alpha u\|_{\mathrm{L}^2(B_{2R}\backslash B_{\frac{R}{2}})}^{\frac{1}{\alpha}}\nonumber\\
&+ R^{4-\frac{12}{p}}\|b\|_{\mathrm{L}^p(B_{2R}\backslash B_{\frac{R}{2}})}^4+ R^{\frac{3(1-\frac{2}{p})(3\beta-1)-3\beta}{2\beta}}\|(u,b)\|_{\mathrm{L}^p(B_{2R}\backslash B_{\frac{R}{2}})}^{3-\frac{1}{\beta}}\|\Lambda^\beta b\|_{\mathrm{L}^2(B_{2R}\backslash B_{\frac{R}{2}})}^{\frac{1}{\beta}}\nonumber\\
&+ R^{\frac{9(1-\frac{2}{p})-5}{2}}\|(u,b)\|_{\mathrm{L}^p(B_{2R}\backslash B_{\frac{R}{2}})}^{3}.
\end{align}
 And  we can notice that $I_3$ and $I_2$ have the same estimation
\begin{align}\label{3.33}
I_3=&2\int_{\mathbb{R}^{3}}(b\cdot u)(b\cdot \nabla\psi_R(x)) dx\leq \frac{C}{R}\|u\|_{\mathrm{L}^3(B_{\frac{3}{2}R}\backslash B_{\frac{3}{4}R})} \|b\|_{\mathrm{L}^3(B_{\frac{3}{2}R}\backslash B_{\frac{3}{4}R})}^2\nonumber\\
&\hbox{(Using the estimation process of $I_{2}$ ......)}\nonumber\\
  \overset{p\geq2}{\lesssim}&R^{\frac{3(1-\frac{2}{p})(2\alpha-1)-2\alpha}{2\alpha}}\|u\|_{\mathrm{L}^p(B_{2R}\backslash B_{\frac{R}{2}})}^{2-\frac{1}{\alpha}}\|\Lambda^\alpha u\|_{\mathrm{L}^2(B_{2R}\backslash B_{\frac{R}{2}})}^{\frac{1}{\alpha}}\nonumber\\
&+ R^{\frac{3(1-\frac{2}{p})(2\beta-1)-\beta}{\beta}}\|b\|_{\mathrm{L}^p(B_{2R}\backslash B_{\frac{R}{2}})}^{4-\frac{2}{\beta}}\|\Lambda^\beta b\|_{\mathrm{L}^2(B_{2R}\backslash B_{\frac{R}{2}})}^{\frac{2}{\beta}}\nonumber\\
&+ R^{\frac{3(1-\frac{2}{p})(2\alpha-1)-4\alpha}{2\alpha}}\|u\|_{\mathrm{L}^p(B_{2R}\backslash B_{\frac{R}{2}})}^{2-\frac{1}{\alpha}}\|\Lambda^\alpha u\|_{\mathrm{L}^2(B_{2R}\backslash B_{\frac{R}{2}})}^{\frac{1}{\alpha}}\nonumber\\
&+ R^{4-\frac{12}{p}}\|b\|_{\mathrm{L}^p(B_{2R}\backslash B_{\frac{R}{2}})}^4+ R^{\frac{3(1-\frac{2}{p})(3\beta-1)-3\beta}{2\beta}}\|(u,b)\|_{\mathrm{L}^p(B_{2R}\backslash B_{\frac{R}{2}})}^{3-\frac{1}{\beta}}\|\Lambda^\beta b\|_{\mathrm{L}^2(B_{2R}\backslash B_{\frac{R}{2}})}^{\frac{1}{\beta}}\nonumber\\
&+ R^{\frac{9(1-\frac{2}{p})-5}{2}}\|(u,b)\|_{\mathrm{L}^p(B_{2R}\backslash B_{\frac{R}{2}})}^{3}
\end{align}
Therefore, we need to take
\begin{equation}\label{e2}
  \left\{
     \begin{array}{ll}
     \frac{3(1-\frac{2}{p})(2\alpha-1)-2\alpha}{2\alpha}\leq 0 , \\
     \frac{3(1-\frac{2}{p})(2\beta-1)-\beta}{\beta}\leq 0 , \\
     \frac{3(1-\frac{2}{p})(2\alpha-1)-4\alpha}{2\alpha}\leq 0 ,\\
     4-\frac{12}{p}\leq 0 ,\\
     \frac{3(1-\frac{2}{p})(3\beta-1)-3\beta}{2\beta}\leq 0 ,\\
     \frac{9(1-\frac{2}{p})-5}{2}\leq 0.
     \end{array}
   \right.
\end{equation}
Thus, $I_2\rightarrow 0$ and $I_3\rightarrow 0$ hold true when $R\rightarrow \infty$.

Furthermore, we apply Eq. \eqref{eq3.5} and the Gagliardo-Nirenberg inequality to estimate $I_4$
\begin{align}\label{3.34}
I_4=&2\int_{\mathbb{R}^{3}}pu\cdot\nabla\psi_R~dx\leq\frac{C}{R}\|u\|_{\mathrm{L}^3(B_{\frac{3}{2}R}\backslash B_{\frac{3}{4}R})}\|p\|_{\mathrm{L}^{\frac{3}{2}}(B_{\frac{3}{2}R}\backslash B_{\frac{3}{4}R})}\nonumber\\
\leq&\frac{C}{R}\|u\|_{\mathrm{L}^3(B_{\frac{3}{2}R}\backslash B_{\frac{3}{4}R})}\left(\|u\|_{\mathrm{L}^3(\mathbb{R}^3)}^2+ \|b\|_{\mathrm{L}^3(\mathbb{R}^3)}^2\right)\nonumber\\
\leq&\frac{C}{R}\|u\|_{\mathrm{L}^2(B_{2R}\backslash B_{\frac{R}{2}})}^{1-\frac{1}{2\alpha}}\|\Lambda^\alpha u\|_{\mathrm{L}^2(B_{2R}\backslash B_{\frac{R}{2}})}^{\frac{1}{2\alpha}}\|u\|_{\mathrm{L}^2(\mathbb{R}^3)}^{2-\frac{1}{\alpha}}\|\Lambda^\alpha u\|_{\mathrm{L}^2(\mathbb{R}^3)}^{\frac{1}{\alpha}}\nonumber\\
&+\frac{C}{R^{\frac{3}{2}}}\|u\|_{\mathrm{L}^2(B_{2R}\backslash B_{\frac{R}{2}})}\|u\|_{\mathrm{L}^2(\mathbb{R}^3)}^{2-\frac{1}{\alpha}}\|\Lambda^\alpha u\|_{\mathrm{L}^2(\mathbb{R}^3)}^{\frac{1}{\alpha}}\nonumber\\
&+\frac{C}{R}\|u\|_{\mathrm{L}^2(B_{2R}\backslash B_{\frac{R}{2}})}^{1-\frac{1}{2\alpha}}\|\Lambda^\alpha u\|_{\mathrm{L}^2(B_{2R}\backslash B_{\frac{R}{2}})}^{\frac{1}{2\alpha}}\|b\|_{\mathrm{L}^2(\mathbb{R}^3)}^{2-\frac{1}{\beta}}\|\Lambda^\beta b\|_{\mathrm{L}^2(\mathbb{R}^3)}^{\frac{1}{\beta}}\nonumber\\
&+\frac{C}{R^{\frac{3}{2}}}\|u\|_{\mathrm{L}^2(B_{2R}\backslash B_{\frac{R}{2}})}\|b\|_{\mathrm{L}^2(\mathbb{R}^3)}^{2-\frac{1}{\beta}}\|\Lambda^\beta b\|_{\mathrm{L}^2(\mathbb{R}^3)}^{\frac{1}{\beta}}\nonumber\\
\leq&\frac{C}{R}\|u\|_{\mathrm{L}^2(B_{2R}\backslash B_{\frac{R}{2}})}^{1-\frac{1}{2\alpha}}\|\Lambda^\alpha u\|_{\mathrm{L}^2(B_{2R}\backslash B_{\frac{R}{2}})}^{\frac{1}{2\alpha}}\|u\|_{\mathrm{L}^2(\mathbb{R}^3)}^{2-\frac{1}{\alpha}}\|\Lambda^\alpha u\|_{\mathrm{L}^2(\mathbb{R}^3)}^{\frac{1}{\alpha}}\nonumber\\
&+\frac{C}{R^{\frac{3}{2}}}\|u\|_{\mathrm{L}^2(B_{2R}\backslash B_{\frac{R}{2}})}\|u\|_{\mathrm{L}^2(\mathbb{R}^3)}^{2-\frac{1}{\alpha}}\|\Lambda^\alpha u\|_{\mathrm{L}^2(\mathbb{R}^3)}^{\frac{1}{\alpha}}\nonumber\\
&+\frac{C}{R}\|u\|_{\mathrm{L}^2(B_{2R}\backslash B_{\frac{R}{2}})}^{2-\frac{1}{\alpha}}+\frac{C}{R}\|\Lambda^\alpha u\|_{\mathrm{L}^2(B_{2R}\backslash B_{\frac{R}{2}})}^{\frac{1}{\alpha}}\|b\|_{\mathrm{L}^2(\mathbb{R}^3)}^{4-\frac{2}{\beta}}\|\Lambda^\beta b\|_{\mathrm{L}^2(\mathbb{R}^3)}^{\frac{2}{\beta}}\nonumber\\
&+\frac{C}{R^{\frac{3}{2}}}\|u\|_{\mathrm{L}^2(B_{2R}\backslash B_{\frac{R}{2}})}\|b\|_{\mathrm{L}^2(\mathbb{R}^3)}^{2-\frac{1}{\beta}}\|\Lambda^\beta b\|_{\mathrm{L}^2(\mathbb{R}^3)}^{\frac{1}{\beta}}\nonumber\\
\overset{(p\geq2)}{\lesssim}&R^{\frac{3(1-\frac{2}{p})(6\alpha-3)}{4\alpha}-1}\|u\|_{\mathrm{L}^p(B_{2R}\backslash B_{\frac{R}{2}})}^{1-\frac{1}{2\alpha}}\|\Lambda^\alpha u\|_{\mathrm{L}^2(B_{2R}\backslash B_{\frac{R}{2}})}^{\frac{1}{2\alpha}}\|u\|_{\mathrm{L}^p(\mathbb{R}^3)}^{2-\frac{1}{\alpha}}\|\Lambda^\alpha u\|_{\mathrm{L}^2(\mathbb{R}^3)}^{\frac{1}{\alpha}}\nonumber\\
&+CR^{\frac{3}{2}[(1-\frac{2}{p})(3-\frac{1}{\alpha})-1]}\|u\|_{\mathrm{L}^p(B_{2R}\backslash B_{\frac{R}{2}})}\|u\|_{\mathrm{L}^p(\mathbb{R}^3)}^{2-\frac{1}{\alpha}}\|\Lambda^\alpha u\|_{\mathrm{L}^2(\mathbb{R}^3)}^{\frac{1}{\alpha}}\nonumber\\
&+ R^{\frac{3(1-\frac{2}{p})(2\beta-1)-\beta}{\beta}}\|b\|_{\mathrm{L}^p(\mathbb{R}^3)}^{4-\frac{2}{\beta}}\|\Lambda^\alpha u\|_{\mathrm{L}^2(B_{2R}\backslash B_{\frac{R}{2}})}^{\frac{1}{\alpha}}\|\Lambda^\beta b\|_{\mathrm{L}^2(\mathbb{R}^3)}^{\frac{2}{\beta}}\nonumber\\
&+CR^{\frac{3}{2}[(1-\frac{2}{p})(3-\frac{1}{\beta})-1]}\|u\|_{\mathrm{L}^p(B_{2R}\backslash B_{\frac{R}{2}})}\|b\|_{\mathrm{L}^p(\mathbb{R}^3)}^{2-\frac{1}{\beta}}\|\Lambda^\beta b\|_{\mathrm{L}^2(\mathbb{R}^3)}^{\frac{1}{\beta}}\nonumber\\
&+ R^{\frac{3(1-\frac{2}{p})(2\alpha-1)-2\alpha}{2\alpha}}\|u\|_{\mathrm{L}^p(B_{2R}\backslash B_{\frac{R}{2}})}^{2-\frac{1}{\alpha}}.
\end{align}
It is obvious that when both constraint conditions \eqref{e1} and \eqref{e2} are satisfied, and add the following condition:
\begin{equation}\label{3.37}
	\left\{
	\begin{array}{ll}
		\frac{3}{2}[(1-\frac{2}{p})(3-\frac{1}{\alpha})-1]\leq 0  & \hbox{$\Leftrightarrow p \leq \frac{6\alpha-2}{2\alpha-1}$,} \\
		\frac{3}{2}[(1-\frac{2}{p})(3-\frac{1}{\beta})-1]\leq 0  & \hbox{$\Leftrightarrow p\leq \frac{6\beta-2}{2\beta-1}$,} \\
	\end{array}
	\right.
\end{equation}
we necessarily have $I_4\rightarrow 0$ as $R\rightarrow \infty$.

Finally, combining the vanishing conditions for $I_1-I_4$, i.e. \eqref{e1},   \eqref{e2} and \eqref{3.19}with \eqref{3.37},  the exponent $p$ must satisfy the following constraint condition:
\begin{equation}\label{3.35}
  \left\{
 \begin{array}{ll}
 	2-\frac{9}{p}\leq 0   &  \hbox{$\Leftrightarrow p \leq \frac{9}{2}$,}\\
 	4-\frac{12}{p}\leq 0  &  \hbox{$\Leftrightarrow p\leq 3$,}\\
 	\frac{3}{2}[(1-\frac{2}{p})(3-\frac{1}{\alpha})-1]\leq 0  & \hbox{$\Leftrightarrow p \leq \frac{6\alpha-2}{2\alpha-1}$,} \\
 	\frac{3}{2}[(1-\frac{2}{p})(3-\frac{1}{\beta})-1]\leq 0  & \hbox{$\Leftrightarrow p\leq \frac{6\beta-2}{2\beta-1}$,} \\
 	\frac{3(1-\frac{2}{p})(3\beta-1)-3\beta}{2\beta}\leq 0  &  \hbox{$\Leftrightarrow p\leq \frac{2(3\beta-1)}{2\beta-1}$,}\\
 	\left.
 	\begin{aligned}
 		\frac{3\left(1-\frac{2}{p}\right)(2\alpha-1)-4\alpha}{2\alpha} &\leq 0 \\
 		\frac{3\left(1-\frac{2}{p}\right)(2\alpha-1)-2\alpha}{2\alpha} &\leq 0
 	\end{aligned}
 	\right\}  &  \hbox{$\Leftrightarrow 2 - \frac{6}{p} + 3(\frac{1}{p} - \frac{1}{2})\frac{1}{\alpha}\leq 0$}\\
 \end{array}
 \right.
\end{equation}
and
\begin{equation}\label{3.36}
     \begin{cases}
      \frac{3(1-\frac{2}{p})(6\alpha-3)}{4\alpha}-1 \leq 0,\\
      \frac{3(1-\frac{2}{p})(2\beta-1)-\beta}{\beta} \leq 0.
 \end{cases}
\end{equation}
For the constraint conditions  \eqref{3.35}, we notice that $\frac{2(3\alpha-1)}{2\alpha-1},~\frac{2(3\beta-1)}{2\beta-1}$ are greater than $3$ when $\alpha,~\beta>\frac{1}{2}$. Consequently,
the constraint conditions \eqref{3.35} can be simplified to $p \leq 3$.

Next, we simplify constraint conditions \eqref{3.36}. According to the range of $\alpha \text{ and } \beta$, the simplification process can be divided into four categories:
\begin{enumerate}[label=(\roman*)]
  \item\label{i} For $\frac{1}{2}\leq\alpha\leq\frac{9}{14}$ and $\frac{1}{2}\leq\beta\leq\frac{3}{5}$, we note that
  $$\frac{3(1-\frac{2}{p})(6\alpha-3)}{4\alpha}-1 \leq 0 \text{~~and~~}  \frac{3(1-\frac{2}{p})(2\beta-1)-\beta}{\beta} \leq 0 \text{ always hold true. }$$
Hence, combining with Eq. \eqref{3.35}, we obtain $2\leq p \leq3$. Then the following cases may be proved in the same way as case \ref{i}, so we omit the details.
  \item For $\frac{1}{2}\leq\alpha\leq\frac{9}{14}$ and $\beta>\frac{3}{5}$,   $ \frac{3(1-\frac{2}{p})(2\beta-1)-\beta}{\beta}\leq0$ implies $ p \leq \frac{6(2\beta-1)}{5\beta-3}$. However, $\frac{6(2\beta-1)}{5\beta-3}>3$ is trivial. Combining  with \eqref{3.35}, we conclude $2 \leq p \leq 3$.
  \item For $\frac{9}{14}<\alpha<\frac{3}{4}$ and $\frac{1}{2}\leq\beta\leq\frac{3}{5}$, $\frac{3(1-\frac{2}{p})(6\alpha-3)}{4\alpha}-1 \leq 0$ implies $p \leq \frac{18(2\alpha-1)}{14\alpha-9}$.
   And  we also have $\frac{18(2\alpha-1)}{14\alpha-9}\geq3$ with $\frac{9}{14}<\alpha<\frac{3}{4}$.   Combining  with \eqref{3.35},  we get $2 \leq p \leq 3$.
  \item For $\frac{9}{14}<\alpha<\frac{3}{4}$ and $\beta>\frac{3}{5}$, we obtain $p\leq\min\{\frac{18(2\alpha-1)}{14\alpha-9}, \frac{6(2\beta-1)}{5\beta-3}\}$. Following the aforementioned procedure, we have $2 \leq p \leq 3$.
\end{enumerate}
Therefore,  the constraint conditions \eqref{3.35} and \eqref{3.36} can be simplified to  $2 \leq p \leq 3$   provided $\frac{1}{2}\leq\alpha<\frac{3}{4}$ and $\frac{1}{2}\leq\beta<1$.

 Considering $\Lambda^\alpha u,~\Lambda^\beta b\in L^2(\mathbb{R}^3)$ and Eq. \eqref{3.19}, when $R\rightarrow\infty$, we can obtain
$$2\bm{C}_{\alpha}\int_{\mathbb{R}^{4}_+}y^{\lambda_\alpha}|\bar{\nabla}u^*|^2~dxdy + 2\bm{C}_{\beta}\int_{\mathbb{R}^{4}_+}y^{\lambda_\beta}|\bar{\nabla}b^*|^2~dxdy\leq0,$$
which implies $\bm{C}_{\alpha}\int_{\mathbb{R}^{4}_+}y^{\lambda_\alpha}|\bar{\nabla}u^*|^2~dxdy=0$ and $\bm{C}_{\beta}\int_{\mathbb{R}^{4}_+}y^{\lambda_\beta}|\bar{\nabla}b^*|^2~dxdy=0$. According to Lemma \ref{le2.3} and $\eqref{Eq2.1}_4$, it is obvious that $u=b=0$ because of $u$ is a smooth function in $L^p(\mathbb{R}^3)$.

In summary, we have $u=b=0$ provided $\frac{1}{2}\leq\alpha,~\beta<1$ with $2\leq p\leq\frac{9}{2}$.

\section{The proof of Theorem \ref{th1.5}}
In the absence of the Hall term, equations \eqref{eq1.3} reduce to the MHD equations \eqref{eq1.1}. Consequently, the proof of Theorem \ref{th1.5} is analogous to that of Theorem \ref{th1.1}.
Moreover, the Hall term satisfies the following identity:
$$\int_{\mathbb{R}^3}\nabla\times((\nabla\times b)\times b)\cdot(\nabla\psi_Rb)=\int_{\mathbb{R}^3}((\nabla\times b)\times b)\cdot(\nabla\psi_R\times b).$$

Since we have the following vector identity:
\begin{equation}\label{eq4.1}
	\int_{\mathbb{R}^3}((\nabla\times b)\times b)\cdot(\nabla\psi_R\times b)=\int_{\mathbb{R}^3}[(\nabla\times b)\cdot\nabla\psi_R](b\cdot b)-[(\nabla\times b)\cdot b](b\cdot\nabla\psi_R)~dx.\nonumber
\end{equation}

We then take the $L^2$ inner products of the first two equations in \eqref{eq1.3} with $\psi_Ru$ and $\psi_Rb$, respectively, sum the resulting identities, and integrate by parts.
In this way, we obtain a local energy inequality for the Hall-MHD equations, which differs from \eqref{eq3.4} by only two additional terms:
\begin{align*}
	&2\bm{C}_{\alpha}\int_{\mathbb{R}^{4}_+}y^{\lambda_\alpha}|\bar{\nabla}u^*|^2~\psi_R(x)\chi_R(y)~dxdy + 2\bm{C}_{\beta}y^{\lambda_\beta}|\bar{\nabla}b^*|^2~\psi_R(x)\chi_R(y)~dxdy\\
	\leq&\int_{\mathbb{R}^{3}}(u\cdot\nabla\psi_{R})|u|^2~dx + \int_{\mathbb{R}^{3}}(u\cdot\nabla\psi_{R})|b|^2~dx+2\int_{\mathbb{R}^{3}}(b\cdot u)(b\cdot \nabla\psi_R(x)) dx\\
	&+2\int_{\mathbb{R}^{3}}pu\cdot\nabla\psi_R~dx+ \bm{C}_{s}\int_{\mathbb{R}^{4}_+}(|u^*|^2 +|b^*|^2)\bar{\nabla}(y^{\alpha}\bar{\nabla}(\psi_R(x)\chi_R(y)))~dx\\
	&+\int_{\mathbb{R}^3}[(\nabla\times b)\cdot\nabla\psi_R](b\cdot b)~dx+\int_{\mathbb{R}^3}[(\nabla\times b)\cdot b](b\cdot\nabla\psi_R)~dx\\
	=&\sum_{i=1}^{6}I_i+\int_{\mathbb{R}^3}[(\nabla\times b)\cdot\nabla\psi_R](b\cdot b)~dx+\int_{\mathbb{R}^3}[(\nabla\times b)\cdot b](b\cdot\nabla\psi_R)~dx\\
	:=&\sum_{i=1}^{8}I_i.
\end{align*}

The estimation of $I_1-I_6$ are the same as those in the proof of Theorem 1.1, so we only need to estimate $I_7$ and $I_8$.
\subsection{$\frac{5}{6}\leq\beta<1$}

To begin with, in order to estimate $I_7$, we apply Gagliardo-Nirenberg inequality to estimate $\|b\|_{\mathrm{L}^\frac{6}{3-2\beta}(B_{\frac{3}{2}R}\backslash B_{\frac{3}{4}R})}$ as follows
\begin{align}\label{eq3.14}
  &\|b\|_{\mathrm{L}^\frac{6}{3-2\beta}(B_{\frac{3}{2}R}\backslash B_{\frac{3}{4}R})}\leq\|b\phi_R\|_{\mathrm{L}^\frac{6}{3-2\beta}(\mathbb{R}^3)}\nonumber\\
\leq& C\|\Lambda^\beta (b\phi_R)\|_{\mathrm{L}^2(\mathbb{R}^3)} \hbox{~~~(by Sobolev embedding)}\nonumber\\
\leq&C\|\Lambda^\beta b\|_{\mathrm{L}^2(B_{2R}\backslash B_{\frac{R}{2}})}\|\phi_R\|_{L^{\infty}(\mathbb{R}^3)}  +  C\|b\|_{L^{2}(B_{2R}\backslash B_{\frac{R}{2}})}\|\Lambda^\beta\phi_R\|_{\mathrm{L}^\infty(\mathbb{R}^3)} \hbox{~(by Lemma \ref{le2.2})}\nonumber\\
\leq&C\|\Lambda^\beta b\|_{\mathrm{L}^2(B_{2R}\backslash B_{\frac{R}{2}})}+\frac{C}{R^\beta}\|b\|_{\mathrm{L}^2(B_{2R}\backslash B_{\frac{R}{2}})}.
\end{align}

Moreover, the procedure of proof carries over to $\|\Lambda^{1-\beta} b\|_{\mathrm{L}^\frac{6}{5-4\beta}(B_{\frac{3}{2}R}\backslash B_{\frac{3}{4}R})}$:
\begin{align}\label{eq4.2}
  &\|\Lambda^{1-\beta} b\|_{\mathrm{L}^\frac{6}{5-4\beta}(B_{\frac{3}{2}R}\backslash B_{\frac{3}{4}R})}\leq\|\Lambda^{1-\beta} b\|_{\mathrm{L}^\frac{6}{5-4\beta}(B_{2R}\backslash B_{\frac{R}{2}})}\nonumber\\
\leq& C\|\Lambda^\beta b\|_{\mathrm{L}^2(B_{2R}\backslash B_{\frac{R}{2}})}\|b\|_{\mathrm{L}^2(B_{2R}\backslash B_{\frac{R}{2}})}^0\hbox{~~~(by Gagliardo-Nirenberg inequality)}\nonumber\\
=&C\|\Lambda^\beta b\|_{\mathrm{L}^2(B_{2R}\backslash B_{\frac{R}{2}})}.
\end{align}
Next, we combine \eqref{eq3.14} and \eqref{eq4.2} and use the H\"{o}lder inequality to estimate $I_7$
\begin{align}\label{eq4.5}
&I_7=\int_{\mathbb{R}^{3}}(\nabla\times b)\cdot\nabla\psi_Rb^2~dx\leq\int_{\mathbb{R}^{3}}\varepsilon_{ijk}\partial_j b_k \partial_i\psi_Rb^2~dx\nonumber\\
\leq&\varepsilon_{ijk}\int_{\mathbb{R}^{3}}\partial_j^{\beta_j} b_k \partial_j^{1-\beta_j}(\partial_i\psi_Rb^2)~dx\leq\int_{\mathbb{R}^{3}}\Lambda^\beta b \Lambda^{1-\beta}(b^2\cdot\nabla\psi_R)~dx\hbox{\quad(by Lemma \ref{le2.6})}\nonumber\\
\leq&C\|\Lambda^\beta b\|_{\mathrm{L}^2(B_{\frac{3}{2}R}\backslash B_{\frac{3}{4}R})}\|\Lambda^{1-\beta}(b^2\cdot \nabla \psi_R)\|_{\mathrm{L}^2(B_{\frac{3}{2}R}\backslash B_{\frac{3}{4}R})}  \nonumber\\
\leq&\|\Lambda^\beta b\|_{\mathrm{L}^2(B_{\frac{3}{2}R}\backslash B_{\frac{3}{4}R})}\left(\|\Lambda^{1-\beta} b^2\|_{\mathrm{L}^{\frac{3}{4-3\beta}}(B_{\frac{3}{2}R}\backslash B_{\frac{3}{4}R})}\|\nabla\psi_R\|_{\mathrm{L}^{\frac{6}{6\beta-5}}(B_{\frac{3}{2}R}\backslash B_{\frac{3}{4}R})}\right)\nonumber\\
&+\|\Lambda^\beta b\|_{\mathrm{L}^2(B_{\frac{3}{2}R}\backslash B_{\frac{3}{4}R})}\left(\|\Lambda^{1-\beta}\nabla\psi_R\|_{\mathrm{L}^{\frac{6}{4\beta-3}}(B_{\frac{3}{2}R}\backslash B_{\frac{3}{4}R})}\|b^2\|_{\mathrm{L}^{\frac{3}{3-2\beta}}(B_{\frac{3}{2}R}\backslash B_{\frac{3}{4}R})}\right)\hbox{\quad(by Lemma \ref{le2.2})}\nonumber\\
\leq&\|\Lambda^\beta b\|_{\mathrm{L}^2(B_{\frac{3}{2}R}\backslash B_{\frac{3}{4}R})}\|b\|_{\mathrm{L}^{\frac{6}{3-2\beta}}(B_{\frac{3}{2}R}\backslash B_{\frac{3}{4}R})}\|\Lambda^{1-\beta}b\|_{\mathrm{L}^{\frac{6}{5-4\beta}}(B_{\frac{3}{2}R}\backslash B_{\frac{3}{4}R})}\|\nabla\psi_R\|_{\mathrm{L}^{\frac{6}{6\beta-5}}(B_{\frac{3}{2}R}\backslash B_{\frac{3}{4}R})}\nonumber\\
&+\|\Lambda^\beta b\|_{\mathrm{L}^2(B_{\frac{3}{2}R}\backslash B_{\frac{3}{4}R})}\|\Lambda^{1-\beta}\nabla\psi_R\|_{\mathrm{L}^{\frac{6}{4\beta-3}}(B_{\frac{3}{2}R}\backslash B_{\frac{3}{4}R})}\|b\|_{\mathrm{L}^{\frac{6}{3-2\beta}}(B_{\frac{3}{2}R}\backslash B_{\frac{3}{4}R})}^2\hbox{\quad(by Lemma \ref{le2.2})}\nonumber\\
\leq&\frac{C}{R^{\frac{7-6\beta}{2}}}\left(\|\Lambda^\beta b\|_{\mathrm{L}^2(B_{2R}\backslash B_{\frac{R}{2}})}+\frac{1}{R^\beta}\|b\|_{\mathrm{L}^2(B_{2R}\backslash B_{\frac{R}{2}})}\right)\|\Lambda^\beta b\|_{\mathrm{L}^2(B_{2R}\backslash B_{\frac{R}{2}})}^2\nonumber\\
&+\frac{C}{R^{\frac{7-6\beta}{2}}}\left(\|\Lambda^\beta b\|_{\mathrm{L}^2(B_{2R}\backslash B_{\frac{R}{2}})}^2+\frac{1}{R^{2\beta}}\|b\|_{\mathrm{L}^2(B_{2R}\backslash B_{\frac{R}{2}})}^2\right)\|\Lambda^\beta b\|_{\mathrm{L}^2(B_{2R}\backslash B_{\frac{R}{2}})}\hbox{~~~(by \eqref{eq3.14} and \eqref{eq4.2})}\nonumber\\
\leq&\frac{C}{R^{\frac{7-6\beta}{2}}}\|\Lambda^\beta b\|_{\mathrm{L}^2(B_{2R}\backslash B_{\frac{R}{2}})}^3+\frac{C}{R^{\frac{7-4\beta}{2}}}\|b\|_{\mathrm{L}^2(B_{2R}\backslash B_{\frac{R}{2}})}\|\Lambda^\beta b\|_{\mathrm{L}^2(B_{2R}\backslash B_{\frac{R}{2}})}^2\nonumber\\
&+\frac{C}{R^{\frac{7-2\beta}{2}}}\|b\|_{\mathrm{L}^2(B_{2R}\backslash B_{\frac{R}{2}})}^2\|\Lambda^\beta b\|_{\mathrm{L}^2(B_{2R}\backslash B_{\frac{R}{2}})}
  \end{align}
for $\frac{5}{6}\leq\beta<1$. Subsequently we apply the H\"{o}lder's inequality in \eqref{eq4.5} to obtain that
  \begin{align}\label{4.6}
  I_7\leq&CR^{\frac{6\beta-7}{2}}\|\Lambda^\beta b\|_{\mathrm{L}^2(B_{2R}\backslash B_{\frac{R}{2}})}^3+CR^{2\beta-\frac{3}{p}-2}\|b\|_{\mathrm{L}^p(B_{2R}\backslash B_{\frac{R}{2}})}\|\Lambda^\beta b\|_{\mathrm{L}^2(B_{2R}\backslash B_{\frac{R}{2}})}^2\nonumber\\
&+CR^{\beta-\frac{6}{p}-\frac{1}{2}}\|b\|_{\mathrm{L}^p(B_{2R}\backslash B_{\frac{R}{2}})}^2\|\Lambda^\beta b\|_{\mathrm{L}^2(B_{2R}\backslash B_{\frac{R}{2}})}
  \end{align}
for $p\geq2$. In the same way, we next estimate $I_8$:
\begin{align}\label{4.7}
	I_8=&\int_{\mathbb{R}^3}[(\nabla\times b)\cdot b](b\cdot\nabla\psi_R)~dx=\int_{\mathbb{R}^{3}}\varepsilon_{ijk}\partial_j b_k b_i b_l\partial_l\psi_R~dx\nonumber\\
	=&\varepsilon_{ijk}\int_{\mathbb{R}^3}\partial_j^\beta b_k\partial_j^{1-\beta} (b_i b_l \partial_l\psi_R)~dx\hbox{\quad(by Lemma \ref{le2.6})}\nonumber\\
	\leq& C\int_{\mathbb{R}^3}\Lambda^\beta b\Lambda^{1-\beta} (b_ib_l \nabla\psi_R)~dx\nonumber\\
	\leq&C\|\Lambda^\beta b\|_{\mathrm{L}^2(B_{\frac{3}{2}R}\backslash B_{\frac{3}{4}R})}\|\Lambda^{1-\beta}(b_ib_l\cdot \nabla \psi_R)\|_{\mathrm{L}^2(B_{\frac{3}{2}R}\backslash B_{\frac{3}{4}R})}  \nonumber\\
	\leq&\|\Lambda^\beta b\|_{\mathrm{L}^2(B_{\frac{3}{2}R}\backslash B_{\frac{3}{4}R})}\left(\|\Lambda^{1-\beta} (b_ib_l)\|_{\mathrm{L}^{\frac{3}{4-3\beta}}(B_{\frac{3}{2}R}\backslash B_{\frac{3}{4}R})}\|\nabla\psi_R\|_{\mathrm{L}^{\frac{6}{6\beta-5}}(B_{\frac{3}{2}R}\backslash B_{\frac{3}{4}R})}\right)\nonumber\\
	&+\|\Lambda^\beta b\|_{\mathrm{L}^2(B_{\frac{3}{2}R}\backslash B_{\frac{3}{4}R})}\left(\|\Lambda^{1-\beta}\nabla\psi_R\|_{\mathrm{L}^{\frac{6}{4\beta-3}}(B_{\frac{3}{2}R}\backslash B_{\frac{3}{4}R})}\|b_ib_l\|_{\mathrm{L}^{\frac{3}{3-2\beta}}(B_{\frac{3}{2}R}\backslash B_{\frac{3}{4}R})}\right)\hbox{\quad(by Lemma \ref{le2.2})}\nonumber\\
	\lesssim&\|\Lambda^\beta b\|_{\mathrm{L}^2(B_{\frac{3}{2}R}\backslash B_{\frac{3}{4}R})}\|b\|_{\mathrm{L}^{\frac{6}{3-2\beta}}(B_{\frac{3}{2}R}\backslash B_{\frac{3}{4}R})}\|\Lambda^{1-\beta}b\|_{\mathrm{L}^{\frac{6}{5-4\beta}}(B_{\frac{3}{2}R}\backslash B_{\frac{3}{4}R})}\|\nabla\psi_R\|_{\mathrm{L}^{\frac{6}{6\beta-5}}(B_{\frac{3}{2}R}\backslash B_{\frac{3}{4}R})}\nonumber\\
	&+\|\Lambda^\beta b\|_{\mathrm{L}^2(B_{\frac{3}{2}R}\backslash B_{\frac{3}{4}R})}\|\Lambda^{1-\beta}\nabla\psi_R\|_{\mathrm{L}^{\frac{6}{4\beta-3}}(B_{\frac{3}{2}R}\backslash B_{\frac{3}{4}R})}\|b\|_{\mathrm{L}^{\frac{6}{3-2\beta}}(B_{\frac{3}{2}R}\backslash B_{\frac{3}{4}R})}^2\hbox{\quad(by Lemma \ref{le2.2})}\nonumber\\
\end{align}
Hence, we can find that the estimation of $I_8$ coincides with that of $I_7$, so we can directly obtain the estimate of $I_8$:
\begin{align}\label{4.8}
	I_8\leq&CR^{\frac{6\beta-7}{2}}\|\Lambda^\beta b\|_{\mathrm{L}^2(B_{2R}\backslash B_{\frac{R}{2}})}^3+CR^{2\beta-\frac{3}{p}-2}\|b\|_{\mathrm{L}^p(B_{2R}\backslash B_{\frac{R}{2}})}\|\Lambda^\beta b\|_{\mathrm{L}^2(B_{2R}\backslash B_{\frac{R}{2}})}^2\nonumber\\
	&+CR^{\beta-\frac{6}{p}-\frac{1}{2}}\|b\|_{\mathrm{L}^p(B_{2R}\backslash B_{\frac{R}{2}})}^2\|\Lambda^\beta b\|_{\mathrm{L}^2(B_{2R}\backslash B_{\frac{R}{2}})}
\end{align}

 Therefore, noticing that $\Lambda^\beta b\in L^2(\mathbb{R}^3)$, when $R\rightarrow\infty$, we need to take
\begin{equation}\label{4.3}
  \left\{
  \begin{array}{ll}
       \frac{6\beta-7}{2}\leq 0  & \hbox{$\Leftrightarrow \beta \leq \frac{7}{6}$,} \\
       2\beta-\frac{3}{p}-2\leq 0  & \hbox{$\Leftrightarrow 2\beta-2\leq \frac{3}{p}$,} \\
       \beta-\frac{6}{p}-\frac{1}{2}\leq 0  & \hbox{$\Leftrightarrow p \leq \frac{12}{2\beta-1}$.} \\
  \end{array}
\right.
\end{equation}
Since $\frac{5}{6}\leq \beta < 1$, it is easy to check that $\frac{6}{3-2\beta}\geq\frac{9}{2}$ and $2\beta-2< \frac{3}{p}$ always hold. Hence, when $\frac{1}{2}\leq \alpha< 1$, with the additional condition from $I_7$ \eqref{4.3} applied to Theorem \ref{th1.1}, by repeating the classification procedure described in Subsection 4.1, we can obtain $2 \leq p \leq \frac{9}{2}$ if we have $\frac{1}{2}\leq\alpha<1$ and $\frac{5}{6}\leq \beta <1$.

Thus, combining the two cases, we can complete the proof of Theorem \ref{th1.5} according to the Lemma \ref{le2.3} and $\eqref{Eq2.1}_4$ as $R\rightarrow \infty$.


\end{document}